\documentclass[preprint,11pt]{elsarticle}

\usepackage[utf8]{inputenc}
\usepackage{amsmath}
\usepackage{amssymb}
\usepackage{amsthm}
\usepackage{mathrsfs}  
\usepackage{graphicx,subcaption}
\usepackage[left=2.2cm,right=2.0cm,top=2.5cm,bottom=2.8cm]{geometry}
\usepackage[table]{xcolor}
\usepackage{boldline}
\usepackage{lineno}
\usepackage{tikz}

\usepackage{stmaryrd}  
\usepackage{fancyhdr}

\newtheorem{remark}{Remark}

\newcommand{\llkh}{ \{\!\!\{ }
\newcommand{\rrkh}{ \}\!\!\} }

\newcommand{\bA}{\textbf{A}}
\newcommand{\bB}{\textbf{B}}
\newcommand{\bQ}{\textbf{Q}}
\newcommand{\bu}{\textbf{u}}
\newcommand{\bF}{\textbf{F}}

\newcommand{\bff}{\textbf{f}}
\newcommand{\bgg}{\textbf{g}}
\newcommand{\bn}{\textbf{n}}
\newcommand{\bxi}{\boldsymbol{\xi}}

\begin{document}

\begin{frontmatter}

\title{A new direct discontinuous Galerkin method with interface correction for two-dimensional compressible Navier-Stokes equations}

\author[label1]{Mustafa E. Danis\fnref{label2}}
\ead{danis@iastate.edu}
\author[label1]{Jue Yan\corref{cor1}\fnref{label3}}
\ead{jyan@iastate.edu}

\cortext[cor1]{Corresponding author}
\fntext[label3]{Research work of the author is supported by National Science Foundation grant DMS-1620335 and Simons Foundation grant 637716.}
\address[label1]{Department of Mathematics, Iowa State University, Ames, 50011, USA}
\address[label2]{Department of Aerospace Engineering, Iowa State University, Ames, 50011, USA}

\begin{abstract}
We propose a new formula for the nonlinear viscous numerical flux and extend the direct discontinuous Galerkin method with interface correction (DDGIC) of Liu and Yan \cite{liuyan2010ddgic} to compressible Navier-Stokes equations. The new DDGIC framework is based on the observation that the nonlinear diffusion can be represented as a sum of multiple individual diffusion processes corresponding to each conserved variable. A set of direction vectors corresponding to each individual diffusion process is defined and approximated by the average value of the numerical solution at the cell interfaces. The new  framework only requires the computation of conserved variables' gradient, which is linear and approximated by the original direct DG numerical flux formula. The proposed method greatly simplifies the implementation, and thus, can be easily extended to general equations and turbulence models.
Numerical experiments with $P_1$, $P_2$, $P_3$ and $P_4$ polynomial approximations are performed to verify the optimal $(k+1)^{th}$ high-order accuracy of the method. The new DDGIC method is shown to be able to accurately calculate physical quantities such as lift, drag, and friction coefficients as well as separation angle and Strouhal number.           
\end{abstract}


\end{frontmatter}


\section{Introduction}\label{Sec_intro}

Solutions of compressible Navier-Stokes (NS) equations may develop complicated structures with small scale features, for which high order approximations are much needed to better resolve the complex flow. Therefore, high order and high resolution methods are often desirable to simulate compressible viscous flows. Some of the well-known high order methods in the literature are finite difference method \cite{Lele1992}, spectral element method \cite{Kopriva-2009}, residual distribution method \cite{Abgrall2017}, spectral volume method \cite{Ekaterinaris2005, Wang2007} and discontinuous Galerkin method \cite{GasserThesis2009, Cockburn2017}. For a comparison among several high order methods, we refer to the review article \cite{Wang-Fidkowski-Abgrall-Bass-Caraeni-Cari-Deconinck-Hartmann2013} and the references therein.
In this article, we present a new direct discontinuous Galerkin method with interface correction to discretize complex nonlinear diffusion terms of compressible NS equations.

Due to its flexibility on geometry, $hp$-adaptivity, ease of obtaining high order approximations, and extremely local data structure for parallel computing, discontinuous Galerkin (DG) finite element method is an attractive choice for compressible NS equations. We refer the reader to review articles and books \cite{Cockburn-book, Shu-DGreview2014, Hesthaven_nodal, Riviere2008, Pietro-Ern2011-DGbook} for detailed background on DG method. For Euler equations, the hyperbolic counterpart of NS equations, DG methods are well developed, see \cite{Cockburn-book, Shu-DGreview2014} as well as the recent high order positivity-preserving limiter of Zhang and Shu \cite{zhang-MPS-Review}. On the other hand, for nonlinear diffusion equations and compressible NS equations, DG methods are less studied due to complex viscous terms. Nonetheless, several different ways of handling the diffusion terms in DG methods have been proposed in the literature. One class of such methods consists of the popular Bassi and Rebay method \cite{BR1, BR4},
the hybridizable DG method \cite{Nguyen-Peraire-Cockburn2011}, the compact DG method \cite{peraire2008compact, Persson2013}, and the recent positivity-preserving DG method \cite{Zhang-2017-NS}. In these methods, either auxiliary variables are introduced and approximated separately or an extra step of lifting operation is involved. Thus, they are relatively expensive. Another class of DG methods includes the interior penalty DG (IPDG) method  \cite{arnold1982interior,wheeler1978elliptic,baker1977finite} and its extensions to compressible NS equations \cite{hartmann2006a, hartmann2006b, hartmann2008}.

Direct discontinuous Galerkin (DDG) method is a \emph{diffusion solver} developed by the author and collaborators in \cite{liuyan2008ddg, liuyan2010ddgic, vidden2013sddg, yan2013new}. The key feature of direct DG method is the introduction of a numerical flux that approximates the solution's gradient across discontinuous element edges. The numerical flux formula involves the solution jump, gradient average and second order derivative jump values. As a matter of fact, direct DG method is closely related to IPDG method and degenerates to IPDG method with piecewise constant and linear polynomial approximations. On the other hand, direct DG method is found to possess a number of advantages over IPDG method with higher order approximations ($k\geq 2$), see \cite{chen2016}, \cite{Zhang-Yan-2017} and \cite{Huang-Yan-2020} for discussions on a maximum principle preserving limiter, superconvergence on approximating solution's gradient, and adapting solution jump and flux jump interface conditions for elliptic interface problems, respectively. Among the four versions of DDG methods \cite{liuyan2008ddg, liuyan2010ddgic, vidden2013sddg, yan2013new}, direct DG method \cite{liuyan2010ddgic} with interface correction (DDGIC)  is the best choice for time dependent parabolic problems whereas symmetric DDG method \cite{vidden2013sddg} is more suitable for elliptic problems. Notice that the first direct DG method \cite{liuyan2008ddg} is not preferred since an order loss is observed on nonuniform meshes, especially with even degree polynomial approximations. 

Previously, DDG methods were implemented for compressible NS equations in a handful of studies. The first discussion was carried out by Kannan and Wang in \cite{KANNAN20102007} in a spectral volume method setting. Recently, in a sequence of articles \cite{yang2016fast, cheng2016DDG, cheng2016DDG_RANS, cheng_liu_liu_luo_2017}, Cheng and Luo et al. applied the first direct DG method \cite{liuyan2008ddg} to compressible NS equations, in which only quadratic polynomials were considered and no interface correction term was added. Regarding the viscous numerical flux, they employed a product-rule approach that is consistent with our method on the continuous level. More recently, Yue et al. \cite{yue_cheng_liu_2017} and Cheng et al. \cite{cheng2018DDGIC} extended the first DDG method to have more interface terms added and developed symmetric and interface correction versions of DDG method for NS equations, respectively. Their method might be thought of as a generalization of the IPDG method of Hartmann and Houston \cite{hartmann2006a} by including jump terms for second order derivatives. In these studies, due to the simplicity of the DDG method, efficient implicit solvers were developed by computing the viscous part of the Jacobian matrix exactly. DDG methods were also employed to solve Reynolds-averaged NS equations for turbulent flows with high Reynolds numbers \cite{cheng2016DDG_RANS, YANG2018216}. In the context of compressible NS equations, DDG methods were further studied in \cite{CHENG2018adaptive, Xiaoquan2019}.

Despite all the efforts mentioned above, the original DDGIC method \cite{liuyan2010ddgic} could not be applied to compressible NS equations directly. This is because an antiderivative of each nonlinear diffusion component is required in the numerical flux formula, and unfortunately, it does not exist for compressible NS equations. 
A similar difficulty arises for the interface correction terms. For general system of conservation equations, it is not clear how to define suitable interface correction terms for the test function. For this reason, we see different definitions for the viscous numerical flux and interface correction terms across the existing DDG methods, each of which is designed to be used in a specific application. 

In this article, we propose a new DDG method with interface correction for compressible NS equations. In the new DDGIC method, antiderivatives of each nonlinear diffusion component are not required, and a suitable algorithm is developed for computing the interface correction terms. Therefore, the new DDGIC method can be conveniently used in different applications. The new idea is based on the following observation: For each equation of the system, we can decompose the complicated multi-variable entangled nonlinear diffusion into multiple individual diffusion processes corresponding to each conserved variable. We then explore the algebraic adjoint property of inner product and introduce a new direction vector on the element edge. Nonlinearity of each diffusion process essentially goes into the corresponding new direction vector defined at the cell interface, which can be approximated by the solution average of the two neighboring cells. To complete the computation of the viscous numerical flux, gradient of the conserved variables are approximated by the numerical flux formula of the original DDG method. Also, the new direction vectors are further applied to define the interface correction terms involving the test function.

One important advantage of the new DDGIC method is that it requires a rather simple implementation. Also, it can be extended to more general equations quite easily. For example, the new DDGIC method can be used as is when different turbulence models are implemented. This might improve the code reuse: The computation of the numerical fluxes of the conserved variables' gradients are independent of the governing equation and the number of equations in that system. It is the definition of new direction vectors that reflects the physics of the underlying governing equation in the numerical solution. Therefore, the new DDGIC method can be seen as a general framework to apply DDGIC method to general conservation laws with nonlinear diffusion.

The paper is organized as follows: In Section \ref{sec:2.1}, we first use a scalar nonlinear diffusion equation to review the original DDGIC method and present the new idea to compute nonlinear diffusion term numerical flux.  In Section \ref{sec:2.2}, we present a detailed discussion on the derivation of the new DDGIC scheme formulation for compressible NS equations. In Section \ref{sec:2.3}, the implementation of boundary conditions are described. In Section \ref{Sec_Num_examples}, the high-order accuracy of the new DDGIC method is demonstrated through several numerical examples. Finally, conclusions are drawn in Section \ref{sec:Conclusion}. 

\section{The new direct discontinuous Galerkin method with interface correction}
In this section, we formally derive the new direct discontinuous Galerkin method with interface correction (DDGIC) for 2-D NS equations. First, we briefly discuss our new DDGIC method for the scalar nonlinear equations. We then explain the derivation of the new DDGIC method for 2-D compressible NS equations in detail.   

\subsection{The new DDGIC method for scalar nonlinear diffusion equations}\label{sec:2.1}
We use the following nonlinear diffusion equation to illustrate the idea of the new DDGIC method that will be applied to NS equations
\begin{equation}\label{Eqn_nonlinear_diff}
	\frac{\partial u}{\partial t}=\nabla\cdot(\bA(u)\nabla u),\hspace{0.5cm} (\textbf{x},t)\in\Omega\times(0,T).
\end{equation} 
Here, the diffusion matrix $\bA(u)$ is assumed to be positive definite. Furthermore, $\Omega\subset\mathbb{R}^2$ represents the computational domain.  
Let $\mathcal{T}_h$ be a shape regular triangular mesh partition of $\Omega$ such that $\overline{\Omega}= \cup_{{K\in \mathcal{T}_h}} K$. We moreover define $h_K$ as the diameter of the circle inscribed in the element $K$. The numerical solution space is defined as 
\begin{equation}\label{Eqn_Scalar_solutionSpace}
	\mathbb{V}_h^k:=\{v_h\in L^2(\Omega): v_h(x,y)\vline_K\in \mathbb{P}_k(K)\}
\end{equation} 
where $\mathbb{P}_k(K)$ is the space of polynomials of degree $k$ in two dimensions. Across the discontinuous element edges, we use following notations to denote the solution jump and average
\begin{equation*}
\llbracket u \rrbracket=u^+-u^-, \hspace{1cm} \llkh u\rrkh=\frac{u^+ + u^-}{2},\quad \forall (x,y) \in \partial K,
\end{equation*}
where $u^+$ and $u^-$ refer to the solution values obtained from the exterior and interior of the current element $K$. In this study, we consider an orthonormal set of basis functions $\{\phi_j\}_{j=1}^{N_b}\in\mathbb{P}_k(K)$ such that
\begin{equation}\label{Eqn_basis_orthogonal}
	\int_K \phi_i(x,y) \phi_j(x,y) \; dxdy =\delta_{ij}, \hspace{0.5cm} \forall K \in \mathcal{T}_h,
\end{equation}
where $\delta_{ij}$ is the Kronecker delta function. The DDGIC solution is defined as $u_h(x,y,t)\in \mathbb{V}_h^k$ such that on cell $K$ we have
\begin{equation}\label{Eqn_sol_exp_scalar}
	u_h(x,y,t)|_K=\sum_{i=1}^{N_b} u_i(t)\phi_i(x,y).
\end{equation}

\vspace{.05in}
\noindent
{\underline{\emph{Original DDGIC method \cite{liuyan2010ddgic} :}}}
\vspace{.15in}

The key component of direct discontinuous Galerkin method \cite{liuyan2008ddg} is the numerical flux concept $ \widehat{\bA(u_h)\nabla u_h}$ introduced to approximate $\bA(u_h)\nabla u_h$ across discontinuous element edges $\partial K$. We multiply (\ref{Eqn_nonlinear_diff}) with test function $v_h\in \mathbb{V}_h^k$, perform integration by parts, add interface correction terms and formally obtain the DDGIC method in \cite{liuyan2010ddgic} as
\begin{equation}\label{Eqn_old_DDGIC}
	\int_{K}(u_h)_t v_h ~dxdy=\int_{\partial K}   \widehat{\bA(u_h)\nabla u_h} \cdot \bn\,v_h \;ds - \int_{K}   \bA(u_h)\nabla u_h \cdot \nabla v_h \;dxdy+ \mbox{InterfaceTerms}, \quad \forall K \in \mathcal{T}_h,
\end{equation}
where $\bn$ is the outward unit normal vector on $\partial K$. In the original DDG method \cite{liuyan2008ddg}, the numerical flux $\widehat{\bA(u_h)\nabla u_h}\cdot\bn$ is defined as 
\begin{equation}\label{Eqn_orig_ddg_flux}
	\widehat{a_{ij}(u_h)(u_h)_{x_j}} = \frac{\beta_0}{h_e} \llbracket b_{ij}(u_h) \rrbracket n_j + \llkh b_{ij}(u_h)_{x_j}\rrkh + \beta_1 h_e \llbracket b_{ij}(u_h)_{x_1 x_j} n_1 + b_{ij}(u_h)_{x_2 x_j} n_2\rrbracket,
\end{equation}
where $a_{ij}(u)$ is the $ij$ component of the diffusion matrix $\bA(u)$, and $b_{ij}(u)$ is an antiderivative of $a_{ij}(u)$ defined as $b_{ij}(u)=\int^u a_{ij}(s) ds$. Furthermore, $h_e$ is the average of the cell diameter $h_K$ and its neighbor cell sharing the same edge $\partial K$.
The interface correction term in (\ref{Eqn_old_DDGIC}) also involves the nonlinear jump terms $\llbracket b_{ij}(u_h)\rrbracket$ and the test function gradient on $\partial K$.  
\begin{remark}
The major drawback of the original DDG \cite{liuyan2008ddg} and DDGIC \cite{liuyan2010ddgic} methods is that $b_{ij}(u)$ might not be calculated explicitly if $a_{ij}(u)$ is a complicated function, which is the case for the energy equation of 2-D compressible NS equations as we will discuss in Section \ref{sec:2.2}. Therefore, the extension of the original DDG method to general systems with nonlinear diffusion has been unclear. Next, we introduce the new idea to address this issue.
\end{remark}

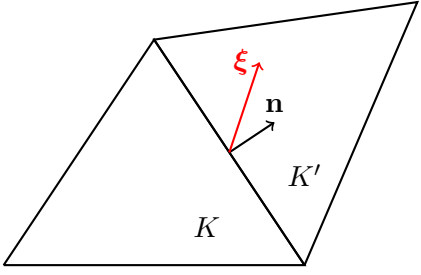
\begin{figure}
    \begin{center}
        \begin{tikzpicture}
            \draw [thick] (0,0) -- (4,0) -- (2,3) -- (0,0);
            \draw [thick] (4,0) -- (2,3) -- (5.5,3.5) -- (4,0);
            \draw [thick,->] (3,1.5) -- (3.6,1.9);
            \draw [thick,red,->] (3,1.5) -- (3.4,2.7);
            \node at (2.7,0.5) {$K$};
            \node at (4,1.2) {$K'$};
            \node at (3.6,1.9) [above] {$\mathbf{n}$};
            \node at (3.4,2.7) [left] {{\color{red}$\bxi$}};
        \end{tikzpicture}
    \caption{New direction vector $\bxi$ on $\partial K$}\label{Fig:new-xi-direction}
    \end{center}
\end{figure}

\vspace{.1in}
\noindent
{\underline{\emph{The new DDGIC method :}}}
\vspace{.1in}

In the new DDGIC method, the same scheme formulation (\ref{Eqn_old_DDGIC})
is used with a new definition for the nonlinear numerical flux $\widehat{\bA(u)\nabla u}\cdot\bn$ and the interface correction term. We take advantage of the adjoint property of the vector inner product 
\begin{equation}\label{Eqn_adjoint_property}
	\bA(u)\nabla u\cdot\bn=\nabla u\cdot \bA(u)^T\bn.
\end{equation}
Then, by rotating and stretching/compressing the unit normal vector $\bn$, we define a new direction vector $\boldsymbol{\xi}=\bA(u)^T\bn$ on edge $\partial K$, see Figure \ref{Fig:new-xi-direction}. Note that this leads to $\bA(u)\nabla u \cdot \bn=\nabla u \cdot \boldsymbol{\xi}(u)=\lVert\bxi\lVert u_\xi$, which is simply proportional to the derivative of $u$ along $\boldsymbol{\xi}$. Here, we denote by $\lVert\:\cdot\:\rVert$ the Euclidean norm. The new numerical flux is defined as
\begin{equation}
\widehat{\bA(u_h)\nabla u_h} \cdot \bn = \widehat{\nabla u_h} \cdot \bxi\left(\llkh u_h \rrkh\right)\quad \mbox{with} \quad \bxi\left(\llkh u_h \rrkh\right)=\bA(\llkh u_h \rrkh)^T\bn, \quad\quad \forall (x,y)\in \partial K,
\end{equation}
and the numerical flux $\widehat{\nabla u_h}=(\widehat{(u_h)_x},\widehat{(u_h)_y})$ is given below as in the original DDGIC method \cite{liuyan2010ddgic}:
\begin{equation}\label{Eqn_scalar_ddg_flux}
\begin{aligned}
\widehat{(u_h)_x}&=\beta_0\frac{\llbracket u_h\rrbracket}{h_e} n_1+\llkh (u_h)_x\rrkh+\beta_1 h_e \llbracket (u_h)_{xx}n_1+(u_h)_{yx}n_2 \rrbracket, \\
\widehat{(u_h)_y}&=\beta_0\frac{\llbracket u_h \rrbracket}{h_e} n_2+\llkh (u_h)_y\rrkh+\beta_1 h_e \llbracket (u_h)_{xy}n_1+(u_h)_{yy}n_2\rrbracket.
\end{aligned}
\end{equation}
A suitable numerical flux coefficient pair $(\beta_0,\beta_1)$ should be chosen to guarantee the stability and optimal convergence of the DDGIC method. In this paper, we simply take $\beta_0=(k+1)^2$ for all $k$, and $\beta_1=\frac{1}{2k(k+1)}$ for $k \ge 2$. 

Over each element $K \in \mathcal{T}_h$, the new DDGIC scheme is then defined as finding the solution $u_h \in \mathbb{V}_h^k$ such that for any test function $v_h \in \mathbb{V}_h^k$ we have
\begin{equation}\label{Eqn_DDGIC_scalar}
\int_{K}(u_h)_t v_h ~dxdy=\int_{\partial K}  \widehat{\nabla u_h} \cdot \bxi(\llkh u_h \rrkh) v_h ~ds - \int_{K} \bA(u_h)\nabla u_h \cdot\nabla v_h ~dxdy 
- \frac{1}{2}\int_{\partial K} \llbracket u_h \rrbracket \nabla v_h \cdot \bxi(\llkh u_h \rrkh) ~ds.
\end{equation} 


In this paper, we adopt orthonormal set of basis functions $\{\phi_j\}_{j=1}^{N_b}$ to write out the DDGIC polynomial solution, see (\ref{Eqn_basis_orthogonal}) and (\ref{Eqn_sol_exp_scalar}). With test function $v_h=\phi_j$ taken in (\ref{Eqn_DDGIC_scalar}), the new DDGIC scheme formulation is simplified as  
\begin{equation*}
	\frac{d{u}_j}{dt}=\int_{\partial K}  \widehat{\nabla u_h} \cdot \bxi(\llkh u_h \rrkh) \,\phi_j ~ds - \int_{K}   \bA(u_h)\nabla u_h \cdot \nabla\phi_j ~dxdy 
- \frac{1}{2}\int_{\partial K} \llbracket u_h \rrbracket \nabla\phi_j \cdot \bxi(\llkh u_h \rrkh) ~ds. 
\end{equation*}
To finish up the full discretization of solving nonlinear diffusion equation with the new DDGIC method, we either employ an explicit method for time discretization, i.e. the strong stability preserving Runge-Kutta method \cite{Gottlieb-Shu-Tabmor-2001,Shu-Osher-1988}, or an  implicit method to solve the above ODE for the solution polynomial coefficient ${u}_j(t)$.

Notice the interface correction term in the original DDGIC method involves the calculation of diffusion matrix with test function, i.e. $\bA=\bA(v_h)$, which is not suitable for 2-D compressible NS equations. Therefore, adjoint property \eqref{Eqn_adjoint_property} is also used to define the interface correction term in (\ref{Eqn_DDGIC_scalar}). Furthermore, if the diffusion process is linear, i.e. $\bA(u)$ is a constant coefficient matrix, the new DDGIC scheme recovers the original DDGIC method \cite{liuyan2010ddgic}. It should also be mentioned that the adjoint property \eqref{Eqn_adjoint_property} was used in the proof of a bound-preserving limiter for the polynomial solution average in \cite{chen2016}. Yet all numerical tests in \cite{chen2016} were carried out with the original DDGIC method (\ref{Eqn_old_DDGIC}) and the numerical flux (\ref{Eqn_orig_ddg_flux}).


\vspace{.1in}
\subsection{The new DDGIC method for 2-D compressible Navier-Stokes equations}\label{sec:2.2}
\vspace{.1in}
The compressible NS equations in an inertial frame with no external force is given as
\begin{equation}\label{Eqn_ns2d}
	\frac{\partial\bQ}{\partial t} + \nabla \cdot \bF_c(\bQ) = \nabla \cdot \bF_v(\bQ,\nabla\bQ), \hspace{0.5cm} (\textbf{x},t)\in\Omega\times(0,T).
\end{equation}
We have $\bQ =\left(\rho, \rho u, \rho v, E\right)^T$ as the vector of conserved variables where $\rho$ is the density, $u$ is the velocity in x-direction, $v$ is the velocity in y-direction and $E$ is the total energy. We denote the velocity vector by $\textbf{u}=(u,v)^T$. Furthermore, we have $\bF_c$ and $\bF_v$ as the convective and viscous fluxes, respectively. The corresponding flux vectors $\bff$ and $\bgg$ in x- and y- directions satisfy the following: 
\begin{equation*}
\nabla \cdot \bF_c(\bQ) =  \frac{\partial \bff_c(\bQ)}{\partial x} + \frac{\partial \bgg_c(\bQ)}{\partial y}, \quad
\nabla \cdot \bF_v(\bQ,\nabla\bQ) = \frac{\partial \bff_v(\bQ,\nabla\bQ)}{\partial x} + \frac{\partial\bgg_v(\bQ,\nabla\bQ)}{\partial y}.
\end{equation*}
The convective fluxes are defined as 
\begin{equation*}
\bff_c(\bQ)=
	\begin{pmatrix}\rho u \\\rho u^2 + p  \\\rho uv \\u(E+p)
	\end{pmatrix},\hspace{1cm}
\bgg_c(\bQ)=
	\begin{pmatrix}\rho u \\ \rho uv  \\ \rho v^2 + p  \\u(E+p)
	\end{pmatrix},
\end{equation*}
where $p$ is the static pressure. Also, we assume a calorically perfect gas and use the ideal gas law as the equation of state to calculate pressure $p=(\gamma-1)\rho e$.
Here, $\gamma$ is the specific heat ratio, which is assumed to be $\gamma=1.4$ throughout this study. In addition, $e$ is the internal energy. It is also related to the total energy as well as the density and the velocity vector by the following identity
\begin{equation*}
	E=\rho e + \frac{1}{2}\rho\left(u^2+v^2\right).
\end{equation*}
Another important variable is the speed of sound, which is given by $a=\sqrt{\frac{\gamma p}{\rho}}$. Furthermore, the viscous fluxes are given by
\begin{equation*}\label{Eqn-NS-viscous}
	\bff_v(\bQ,\nabla\bQ)=
		\begin{pmatrix} 0 \\ \tau_{11}\\ \tau_{21}\\
			u\tau_{11}+v\tau_{21}+\frac{\mu\gamma}{Pr}e_x
		\end{pmatrix},
	\hspace{1cm}
	\bgg_v(\bQ,\nabla\bQ)=
		\begin{pmatrix} 	0 \\ \tau_{12} \\ \tau_{22}  \\
			u\tau_{12}+v\tau_{22}+\frac{\mu\gamma}{Pr}e_y
		\end{pmatrix}.
\end{equation*}
By subscripts $x$ or $y$, we refer to the partial derivative in that respective direction. Moreover, the viscous stress tensor  $\boldsymbol{\tau}$ is defined as
\begin{equation*}
	\boldsymbol{\tau} = 
		\mu
		\begin{pmatrix}
			\frac{4}{3}u_x - \frac{2}{3}v_y &  u_y + v_x \\
			u_y + v_x &  -\frac{2}{3}u_x + \frac{4}{3}v_y \\
		\end{pmatrix}, 
\end{equation*}
where $\mu$ is the dynamic viscosity, and we denote by $\tau_{ij}$ for $i,j=1,2$ the components of the stress tensor $\boldsymbol{\tau}$. Unless otherwise stated, the dynamic viscosity is obtained by the Sutherland's formula $\mu = \mu_{ref} \left(\frac{T}{T_{ref}}\right)^{3/2} \frac{T_{ref}+C_s}{T+C_s}$, 
where $T$ is the temperature, $T_{ref}$ is the reference temperature, $C_s$ is the Sutherland's temperature, and $\mu_{ref}$ is the reference viscosity at $T_{ref}$. Note that temperature $T$ is related to the internal energy by the thermodynamical identity $e=c_vT$, 
where $c_v$ is the specific heat for a constant volume.
Finally, $Pr$ is the nondimensional Prandtl number, which is typically taken to be 0.72 for air. Along with these definitions, the complete set of governing equations is now closed.

In order to derive the DDG scheme formulation for 2-D compressible NS equations, we employ the same triangulation $\mathcal{T}_h$, numerical solution and test function space $\mathbb{V}_h^k$ in (\ref{Eqn_Scalar_solutionSpace}) and set of orthonormal basis functions $\{\phi_i\}_{i=1}^{N_b}$ in \eqref{Eqn_basis_orthogonal} as in Section \ref{sec:2.1}. 

Similar to the numerical solution $u_h(x,y,t)$ in \eqref{Eqn_sol_exp_scalar} for the scalar case, 
we denote by $\bQ_{h}(x,y,t)\in \left(\mathbb{V}_h^k\right)^4$ the numerical solution to 2-D compressible NS equations (\ref{Eqn_ns2d}). The polynomial expansion of the piecewise numerical solution vector on the orthonormal basis $\{\phi_i\}_{i=1}^{N_b}$ is then given as
\begin{equation*}
	\bQ_{h}(x,y,t)|_{K}=\sum_{i=1}^{N_b} \bQ_{i}(t)\phi_i(x,y).
\end{equation*} 

As in \eqref{Eqn_ns2d}, the components of the conserved variable vector have the same meanings, i.e. $\bQ_h^{(1)}=\rho_h$, $\bQ_h^{(2)}=(\rho u)_h$, $\bQ_h^{(3)}=(\rho v)_h$, and $\bQ_h^{(4)}=E_h$. In a word, $\bQ_{h}(x,y,t)|_{K}$ is a four-component piecewise polynomial vector that is expanded on the orthonormal basis $\{\phi_i\}_{i=1}^{N_b}$ of the cell $K$. Note also that the vector of polynomial coefficients $\bQ_{i}(t)$ is time-dependent. 

In order to simplify the presentation of the new direct DG method, we skip the weak formulation of (\ref{Eqn_ns2d}) with a general test function. Instead, we multiply each equation of the system (\ref{Eqn_ns2d}) with a test function $\phi_j$ from the orthonormal basis, integrate over any element $K\in\mathcal{T}_h$, perform the integration by parts and formally obtain the time-evolution equation of the polynomial coefficients $\bQ_{i}(t)$ as 
\begin{equation}\label{Eqn_weak_ns2d}
\frac{d\bQ_{j}}{dt}=\int_K \left(\bF_c(\bQ_{h}) - \bF_v(\bQ_{h},\nabla\bQ_{h})\right) \cdot \nabla\phi_j(x,y)~dxdy 
-\int_{\partial K}\left(\widehat{\bF_c}(\bQ_{h})-\widehat{\bF_v}(\bQ_{h},\nabla\bQ_{h})\right)\cdot\bn \phi_j(x,y)~ds.
\end{equation}
Here $\widehat{\bF_c}$ and $\widehat{\bF_v}$ are the convective and viscous numerical fluxes defined suitably on cell interfaces $\partial K$. Notice that both $\bF_c, \bF_v $ and $\widehat{\bF_c}, \widehat{\bF_v}$ in  (\ref{Eqn_weak_ns2d}) are $4\times 2$ matrix functions and the inner product is interpreted as being carried out between the row vector of the matrix to $\nabla\phi_j$ or the outward pointing normal vector $\bn$ on $\partial K$. The interface correction term is dropped out of the right-hand side of (\ref{Eqn_weak_ns2d}) to simplify the derivation of the new direct DG method. 

As a discontinuous Galerkin method, the numerical solution is discontinuous at the cell interfaces. It is the numerical flux that determines how a cell exchanges information with its neighboring cells. Among the many methods in the literature to compute the convective numerical flux, we choose the local Lax-Friedrichs method, which is written out as
\begin{equation}\label{eq:LaxFriedrichs}
\widehat{\bF_c}\left(\bQ_h\right)\cdot\bn=\llkh\bff_c\left(\bQ_h\right)\rrkh n_1 + \llkh\bgg_c\left(\bQ_h\right)\rrkh n_2 - \alpha \llbracket \bQ_h \rrbracket,\quad \alpha= \max\left(\lVert\bu\rVert+a\right).
\end{equation}
where $a$ is the speed of sound.

In the remainder of this section, we will focus on the development of the viscous numerical flux $\widehat{\bF_v}$ and the interface correction term of (\ref{Eqn_weak_ns2d}). First, we will consider the application of adjoint property to the system case on the continuous level, which is a critical step towards the new scheme formulation. We will then present the new DDGIC method for 2-D compressible NS equations (\ref{Eqn_ns2d}) along with the boundary condition implementation. 
\subsubsection{Adjoint property for the multi-variable system case of compressible Navier-Stokes}
\vspace{.1in}
As discussed in section {\ref{sec:2.1}}, the idea of the new DDGIC method is based on the adjoint property 
$$
\bA(u)\nabla u\cdot\bn=\nabla u\cdot \bA(u)^T\bn,
$$
of the nonlinear diffusion (\ref{Eqn_nonlinear_diff}) matrix and the corresponding vector inner product. For the system case of compressible NS equation (\ref{Eqn_ns2d}), we first rewrite each row vector of $(\bff_v(\bQ,\nabla\bQ),\bgg_v(\bQ,\nabla\bQ))$ as the summation of the product of individual nonlinear diffusion matrix and the gradient of its state variable. Then, the new DDGIC scheme for 2-D compressible NS equations will be a straightforward extension of (\ref{Eqn_DDGIC_scalar}). Let's use the x-momentum equation to illustrate the symbolic deduction. The second row of $(\bff_v(\bQ,\nabla\bQ),\bgg_v(\bQ,\nabla\bQ))$ consists of the stress tensor components $(\tau_{11},\tau_{12})$ that can be decomposed into and laid out in column vector format as
\begin{equation*}
	\begin{aligned}
	\begin{pmatrix}	\tau_{11}  \\ \tau_{12}	\end{pmatrix} =
		\mu	
		\begin{pmatrix}	\frac{4}{3}u_x - \frac{2}{3}v_y  \\	u_y + v_x	
		\end{pmatrix}
		&= \mu
		\begin{pmatrix}\frac{4}{3}u_x  \\	u_y
		\end{pmatrix}
		+\mu
		\begin{pmatrix}	- \frac{2}{3}v_y  \\v_x
		\end{pmatrix} 
		=  \mu
		\begin{pmatrix} \frac{4}{3} & 0 \\	0 & 1
		\end{pmatrix} \nabla u
		+\mu
		\begin{pmatrix}0 & -\frac{2}{3}  \\1 & 0
		\end{pmatrix} \nabla v \\
		&=	-\frac{\mu}{\rho}
		\begin{pmatrix}\frac{4}{3}u & -\frac{2}{3}v \\v & u
		\end{pmatrix} \nabla\rho +
		\frac{\mu}{\rho}
		\begin{pmatrix}\frac{4}{3} & 0 \\ 0 & 1
		\end{pmatrix} \nabla (\rho u)  +
		\frac{\mu}{\rho}
		\begin{pmatrix}0 & -\frac{2}{3} \\1 & 0
		\end{pmatrix} \nabla (\rho v).
	\end{aligned}
\end{equation*}
We use the standard ordering and designate the x-momentum equation as the \textit{second} equation, and adopt $\rho=\bQ^{(1)}$, $\rho u=\bQ^{(2)}$, $\rho v=\bQ^{(3)}$, and $E=\bQ^{(4)}$ as the \textit{first}, \textit{second}, \textit{third} and \textit{fourth} conserved variables. For the x-momentum equation, we introduce the following nonlinear diffusion matrices: 
\begin{equation}\label{Eqn_x_mom_A-matrix}
	\bA^{(21)}(\bQ) = -\frac{\mu}{\rho}
	\begin{pmatrix}	\frac{4}{3}u & -\frac{2}{3}v \\ v & u\end{pmatrix}, \hspace{0.4cm}
	\bA^{(22)}(\bQ) =	\frac{\mu}{\rho}
	\begin{pmatrix} \frac{4}{3} & 0 \\  0 & 1\end{pmatrix},\hspace{0.4cm}
	\bA^{(23)}(\bQ)=\frac{\mu}{\rho}
	 \begin{pmatrix}  0 & -\frac{2}{3} \\ 1 & 0	\end{pmatrix}, 
\end{equation}
and the corresponding viscous term can be rewritten as
\begin{equation}\label{Eqn_x_mom_rewritten}
   \begin{pmatrix}	\tau_{11}  \\ \tau_{12}	\end{pmatrix} =
	\mu
	\begin{pmatrix}
		\frac{4}{3}u_x - \frac{2}{3}v_y  \\
		u_y + v_x
	\end{pmatrix}
	= \bA^{(21)}(\bQ)\nabla\rho + \bA^{(22)}(\bQ)\nabla(\rho u) + \bA^{(23)}(\bQ)\nabla(\rho v).
\end{equation}
Similar to the scalar nonlinear diffusion equation, we explore the adjoint property of vector inner product 
$$\bA^{(2m)}(\bQ)\nabla\bQ^{(m)}\cdot\bn=\nabla\bQ^{(m)}\cdot \bA^{(2m)}(\bQ)^T\bn,\quad \mbox{for}\,\, m=1, 2, 3.$$
Correspondingly, the following new direction vectors are defined
\begin{equation}\label{Eqn_x_mom_xi}
	\bxi^{(2m)}(\bQ)=\bA^{(2m)}(\bQ)^T \bn,
\end{equation}
and these lead to
\begin{equation}\label{Eqn_x_mom_adjoint}
	\bA^{(2m)}(\bQ)\nabla\bQ^{(m)}\cdot\bn =\nabla\bQ^{(m)}\cdot\bxi^{(2m)}(\bQ)\quad \mbox{for}\,\, m=1, 2, 3.
\end{equation}
Since there are no viscous terms for the continuity equation, we have $\bA^{(1m)}(\bQ)=\bf 0$ with $m=1,\cdots,4$. For the y-momentum and the total energy equations, we refer to \ref{sec:appendix-A} where a complete description of the diffusion matrices $\bA^{(3m)}(\bQ)$ and $\bA^{(4m)}(\bQ)$ for $m=1,\cdots,4$ is given, respectively.

\begin{remark}
The viscous numerical flux $\widehat{\bF_v}\cdot\bn$ can be computed by the original DDG numerical flux formula \eqref{Eqn_orig_ddg_flux} for each $\widehat{\bA^{(lm)}(\bQ_h)\nabla\bQ_h^{(m)}}$ term if and only if the compatibility condition
\begin{equation}
	\frac{\partial\bB^{(l)}(\bQ)}{\partial\bQ^{(m)}}=\bA^{(lm)}(\bQ)\quad \mbox{for}\,\, m=1,\cdots,4
\end{equation}
is satisfied. Here, $\bB^{(l)}(\bQ)$ is the antiderivative matrix corresponding to the diffusion matrices $\bA^{(lm)}(\bQ)$.
\end{remark}

\begin{remark}
For 2-D compressible NS equations, the diffusion matrices corresponding to x- and y-momentum equations of 2-D compressible NS equations are compatible, i.e. $\bB^{(2)}(\bQ)$ and $\bB^{(3)}(\bQ)$ exist. However, the matrix antiderivative $\bB^{(4)}(\bQ)$ does not exist for the energy equation. Therefore, the energy equation is not compatible and the original DDG method \cite{liuyan2008ddg} cannot be used for 2-D compressible NS equation. Thus, another route will have to be taken.
\end{remark}

\subsubsection{The new DDGIC scheme formulation for 2-D compressible Navier-Stokes}

\vspace{.1in}

We now present the derivation of the new DDGIC method for 2-D NS equations. For the sake of simplicity, we follow an equation-by-equation approach. Since there are no diffusion terms in the continuity equation, we begin by the x-momentum equation. The main idea of the new direct DG method is based on the adjoint property of vector inner product. Motivated by the continuous level adjoint property (\ref{Eqn_x_mom_adjoint}) with the new direction vector $\bxi^{(2m)}(\bQ)$ in (\ref{Eqn_x_mom_xi}) and the corresponding diffusion matrices $\bA^{(2m)}(\bQ)$ in (\ref{Eqn_x_mom_A-matrix}), the viscous numerical flux for x-momentum equation in \eqref{Eqn_weak_ns2d} is defined as
\begin{equation}\label{Eqn_x_mom_diff_numFlux}
\widehat{\bF^{(2)}_v}\cdot\bn= 
\mu	\widehat{
	\begin{pmatrix}
	\frac{4}{3}u_x - \frac{2}{3}v_y  \\	u_y + v_x
	\end{pmatrix}}\cdot\bn	=
	\widehat{\nabla\rho_h}\cdot\bxi^{(21)}(\llkh \bQ_h \rrkh) + \widehat{\nabla(\rho u)_h} \cdot\bxi^{(22)}(\llkh \bQ_h \rrkh)+ \widehat{\nabla(\rho v)_h}\cdot\bxi^{(23)}(\llkh \bQ_h \rrkh).
\end{equation}
Following the adjoint property, the new direction vectors are defined as below 
\begin{equation}
\bxi^{(2m)}(\llkh \bQ_h\rrkh))=\bA^{(2m)}(\llkh \bQ_h\rrkh)^T\bn, \quad m=1, 2, 3,\quad \forall \partial K.
\end{equation}
Notice the notation $\llkh \bQ_h\rrkh=\frac{\bQ^-_h+\bQ^+_h}{2}$ denotes the solution average on $\partial K$ with $\bQ_{h}^-$ and $\bQ_{h}^+$ representing the values of the solution polynomial evaluated from the interior and exterior of the current cell $K$. The simple numerical flux formula of the direct DG method is then employed to compute $\widehat{\nabla \rho_h}=\widehat{\nabla\bQ_h^{(1)}}$, $\widehat{\nabla (\rho u)_h}=\widehat{\nabla\bQ_h^{(2)}}$, $\widehat{\nabla (\rho v)_h}=\widehat{\nabla\bQ_h^{(3)}}$ and $\widehat{\nabla E_h}=\widehat{\nabla\bQ_h^{(4)}}$ on $\partial K$ as follows 
\begin{equation}\label{eq:numvisc_flux_bQ}
	\begin{aligned}
		\widehat{\partial_x\bQ^{(m)}_h} &= \frac{\beta_0}{h_e} \llbracket\bQ_h^{(m)}\rrbracket n_1+\llkh\partial_x\bQ_h^{(m)}\rrkh+\beta_1 h_e \llbracket\partial_x\partial_x\bQ_h^{(m)}n_1+\partial_y\partial_x\bQ_h^{(m)}n_2\rrbracket \\
		\widehat{\partial_y\bQ^{(m)}_h} &= \frac{\beta_0}{h_e}\llbracket\bQ_h^{(m)}\rrbracket n_2+\llkh\partial_y\bQ_h^{(m)}\rrkh+\beta_1 h_e \llbracket\partial_x\partial_y\bQ_h^{(m)}n_1+\partial_y\partial_y\bQ_h^{(m)}n_2\rrbracket
	\end{aligned},\quad m=1, 2, 3, 4,
\end{equation}
which are the vector counterparts of the numerical flux (\ref{Eqn_scalar_ddg_flux}) in the scalar case. Repeating the same process for y-momentum and total energy equations shows that this structure is maintained. 
That is, the following new direction vectors are introduced on $\partial K$
\begin{equation}\label{eq:adj_prop_ns}
    \bxi^{(lm)}=\bxi^{(lm)}(\llkh \bQ_h \rrkh)=\bA^{(lm)}(\llkh \bQ_h \rrkh)^T\bn,\quad  l=2, 3, 4,\quad m=1,\cdots, 4,
\end{equation} 
and then employed along with the numerical flux formula \eqref{eq:numvisc_flux_bQ} to complete the computation of the viscous numerical flux $\widehat{\bF_v}\cdot\bn$. Again, the reader is referred to \ref{sec:appendix-A} for a full description of the diffusion matrices $\bA^{(lm)}$.

It should be noted that the adjoint property \eqref{eq:adj_prop_ns} allows the essentially difficult nonlinear part of the numerical flux to be dealt with the direction vectors $\bxi^{(lm)}$. Then, it is combined with the simple numerical flux formula of the direct DG method to approximate the gradient of conserved variables. In summary, the new DDGIC method defines multiple individual diffusion processes, each of which is combined to calculate the viscous numerical fluxes for 2-D compressible NS equations as 
\begin{equation}\label{eq:Fvhat_NS}
	\widehat{\bF_v}(\bQ_{h},\nabla\bQ_{h})\cdot\bn=
	\begin{pmatrix}
		0\\
		\widehat{\nabla\rho_h}\cdot \bxi^{(21)}+\widehat{\nabla(\rho u)_h}\cdot \bxi^{(22)}+\widehat{\nabla(\rho v)_h}\cdot \bxi^{(23)}\\
		\widehat{\nabla\rho_h}\cdot \bxi^{(31)}+\widehat{\nabla(\rho u)_h}\cdot \bxi^{(32)}+\widehat{\nabla(\rho v)_h}\cdot \bxi^{(33)}\\
		\widehat{\nabla\rho_h}\cdot \bxi^{(41)}+\widehat{\nabla(\rho u)_h}\cdot \bxi^{(42)}+\widehat{\nabla(\rho v)_h}\cdot \bxi^{(43)}+\widehat{\nabla E_h} \cdot \bxi^{(44)}\\
	\end{pmatrix}.
\end{equation}
Similarly, the interface correction term can be simplified and computed as follows
\begin{equation}\label{eq:Icorr_NS}
	\textbf{I}_{corr}(\bQ_h,\phi_j)\cdot\bn=\frac{1}{2}
	\begin{pmatrix}
		0 \\
		\left(\llbracket\rho_h\rrbracket \bxi^{(21)}+ \llbracket (\rho u)_h\rrbracket  \bxi^{(22)} + \llbracket (\rho v)_h\rrbracket  \bxi^{(23)}\right)\cdot\nabla\phi_j(x,y) \\
		\left(\llbracket \rho_h\rrbracket  \bxi^{(31)} + \llbracket (\rho u)_h\rrbracket  \bxi^{(32)} + \llbracket (\rho v)_h\rrbracket  \bxi^{(33)}\right)\cdot\nabla\phi_j(x,y) \\
		\left(\llbracket \rho_h\rrbracket  \bxi^{(41)} + \llbracket (\rho u)_h\rrbracket  \bxi^{(42)} + \llbracket (\rho v)_h\rrbracket  \bxi^{(43)} + \llbracket E_h\rrbracket  \bxi^{(44)}\right)\cdot\nabla\phi_j(x,y)
	\end{pmatrix}.
\end{equation} 
Finally, we write the new DDGIC scheme for 2-D compressible NS equations

\begin{equation}\label{Eqn_DDGIC_ns2d}
	\begin{aligned}
		\frac{d\hat{\bQ}_{j}}{dt}&=\int_K \left(\bF_c(\bQ_{h}) - \bF_v(\bQ_{h},\nabla\bQ_{h})\right) \cdot \nabla\phi_j(x,y)\;dxdy      \\ 
&-\int_{\partial K} \left(\widehat{\bF_c}(\bQ_{h})
-\widehat{\bF_v}(\bQ_{h},\nabla\bQ_{h})\right)\cdot\phi_j(x,y)\bn ~ds
- \int_{\partial K} \textbf{I}_{corr}(\bQ_h,\phi_j)\cdot\bn ~ds, \,\,\forall K \in T_h  
	\end{aligned}
\end{equation}	

One advantage of the proposed approach is that its implementation is very convenient and straightforward. It allows reuse of the code when simulating different equations, i.e turbulence simulations. The new DDG method only requires the computation of the numerical fluxes of the conserved variables' gradients, which are completely independent of the governing equation. The information pertaining to the governing equation enters the numerical flux through the direction vectors $\bxi^{(lm)}$. Therefore, the new DDGIC method puts forth a general framework for applying DDGIC method for general conservation laws with nonlinear diffusion.

\begin{remark}
For a general system of conservation laws of $N$ equations with nonlinear diffusion, viscous numerical flux can be approximated as 
\begin{equation*}
	\widehat{\bF_v}(\bQ_{h},\nabla\bQ_{h})\cdot\bn=
	\sum_{n=1}^{N}
	\begin{pmatrix}
		\widehat{\nabla\bQ^{(n)}_h}\cdot\bxi^{(1n)}(\llkh \bQ_h \rrkh) \\
		\widehat{\nabla\bQ^{(n)}_h}\cdot\bxi^{(2n)}(\llkh \bQ_h \rrkh) \\
		\vdots \\
		\widehat{\nabla\bQ^{(n)}_h}\cdot\bxi^{(Nn)}(\llkh \bQ_h \rrkh)
	\end{pmatrix},
\end{equation*}    
while the interface correction term is given by
\begin{equation*}
	\textbf{I}_{corr}(\bQ_h,\phi_j)\cdot\bn=
	\frac{1}{2}
	\sum_{n=1}^{N}
		\llbracket \nabla\bQ^{(n)}_h\rrbracket 
		\begin{pmatrix}
			\bxi^{(1n)}(\llkh \bQ_h \rrkh)\cdot\nabla\phi_j \\
			\bxi^{(2n)}(\llkh \bQ_h \rrkh)\cdot\nabla\phi_j \\
			\vdots                                  \\
			\bxi^{(Nn)}(\llkh \bQ_h \rrkh)\cdot\nabla\phi_j
	\end{pmatrix}
\end{equation*}    
\end{remark}

\subsection{Boundary Conditions} \label{sec:2.3}
In this section, we describe the boundary condition implementation with the new DDGIC method. All boundary conditions are enforced weakly through the numerical fluxes of  \eqref{eq:LaxFriedrichs} and \eqref{eq:Fvhat_NS} and the interface correction \eqref{eq:Icorr_NS}. 

For any edge $\partial K$ falling on the domain boundary $\partial\Omega$, we denote the interior state by the superscript $-$ and the ghost state by $+$. That is, $\bQ_h^-,\partial_{x_i}\bQ_h^-$ and $\partial_{x_i}\partial_{x_j}\bQ_h^-$ are calculated from the corresponding interior element using the solution polynomial. On the other hand, $\bQ_h^+,\partial_{x_i}\bQ_h^+$ and $\partial_{x_i}\partial_{x_j}\bQ_h^+$ are computed based on the particular boundary condition type, which are discussed below. 

\subsubsection{Periodic boundary}
We simply take the solution and derivatives from the interior state of the corresponding periodic boundary. For instance, if the solution is periodic at $(x,y)$ with $\bQ(x,y)=\bQ(x+L_x,y+L_y)$, then we set
\begin{equation}\label{BC:periodic}
    \begin{aligned}
        \bQ_h^+(x,y) &= \bQ_h^-(x+L_x,y+L_y),\\
        \partial_{x_i}\bQ_h^+(x,y) &=\partial_{x_i}\bQ_h^-(x+L_x,y+L_y), \\
        \partial_{x_i}\partial_{x_j}\bQ_h^+(x,y) &=\partial_{x_i}\partial_{x_j}\bQ_h^- (x+L_x,y+L_y).
    \end{aligned}
\end{equation}

\subsubsection{Inflow and farfield boundaries}
$\bQ_h^+$ is determined from the free-stream conditions while $\partial_{x_i}\bQ_h^+$ and $\partial_{x_i}\partial_{x_j}\bQ_h^+$ are set equal to the corresponding interior values.
\begin{equation}\label{BC:inflow_farfield}
        \bQ^+_h =
        \begin{pmatrix}
            \rho_\infty \\
            \rho_\infty u_\infty \\
            \rho_\infty v_\infty \\
             \frac{p_\infty}{\gamma - 1} + \frac{1}{2}\rho_\infty(u_\infty^2 + v_\infty^2) \\
        \end{pmatrix},\quad 
        \partial_{x_i}\bQ_h^+ =\partial_{x_i}\bQ_h^-,\quad 
        \partial_{x_i}\partial_{x_j}\bQ_h^+ =\partial_{x_i}\partial_{x_j}\bQ_h^-.
\end{equation}
\begin{remark}
Note that this boundary condition is well-posed since the characteristic wave directions are taken into account by the local Lax-Friedrichs flux
\eqref{eq:LaxFriedrichs}. We refer to \cite{nektar_code} for more details.
\end{remark}

\begin{remark}
Since the first and second derivatives of the ghost state in (\ref{BC:inflow_farfield}) are the same as those of the interior state, the numerical flux \eqref{eq:numvisc_flux_bQ} degenerates into
\begin{equation}
	\begin{aligned}
		\widehat{\partial_x\bQ_h} &= \frac{\beta_0}{h_e} \llbracket\bQ_h\rrbracket n_1+\partial_x\bQ_h^-, \\
		\widehat{\partial_y\bQ_h} &= \frac{\beta_0}{h_e}\llbracket\bQ_h\rrbracket n_2+\partial_y\bQ_h^-.
	\end{aligned}
\end{equation}
\end{remark}

\subsubsection{Outflow boundary}
$\bQ_h^+$ is determined by the partially non-reflecting pressure outflow (PNR) method described in \cite{nektar_code} while $\partial_{x_i}\bQ_h^+$ and $\partial_{x_i}\partial_{x_j}\bQ_h^+$ are set equal to the corresponding interior values.
\begin{equation}\label{BC:outflow}
        \bQ^+_h =
        \begin{pmatrix}
            \rho_h^- \\
            (\rho u)_h^- \\
            (\rho v)_h^- \\
             \frac{2p_h^+-p_h^-}{\gamma - 1} + \frac{1}{2}\rho_h^-\left((u_h^-)^2 + (v_h^-)^2\right) \\
        \end{pmatrix}, \quad
        \partial_{x_i}\bQ_h^+ =\partial_{x_i}\bQ_h^-,\quad
        \partial_{x_i}\partial_{x_j}\bQ_h^+ =\partial_{x_i}\partial_{x_j}\bQ_h^-.
\end{equation}

\subsubsection{Adiabatic viscous wall}
$\bQ_h^+$ is computed by assuming no-slip velocity  $(u_h^+,v_h^+)^T=(0,0)^T$ and all other variables are the same as the interior state. On the other hand, $\partial_{x_i}\bQ_h^+$ is computed by assuming zero heat flux at the wall. For a calorically perfect gas, the wall heat flux is zero if the normal derivative of internal energy is zero, i.e. $\partial_n e_h^+=0$. This is achieved by taking projection of the internal energy gradient on the wall, i.e. $\partial_{x_i} e_h^+ = \partial_{x_i} e_h^- - \partial_n e_h^- n_i$. All other variables are the same as the interior state. Moreover, $\partial_{x_i}\partial_{x_j}\bQ_h^+$ is set equal to the interior second derivative.

\begin{equation}\label{BC:adiabatic_wall}
	\begin{aligned}
        \bQ^+_h &=
        \begin{pmatrix}
            \rho_h^- \\
            0 \\
            0 \\
             \frac{2p_h^+-p_h^-}{\gamma - 1} + \frac{1}{2}\rho_h^-\left((u_h^-)^2 + (v_h^-)^2\right) \\
        \end{pmatrix},\\
        \partial_{x_i}\bQ_h^+ &=
        \begin{pmatrix}
            \partial_{x_i} \rho_h^- \\
            \partial_{x_i} (\rho u)_h^- \\
            \partial_{x_i} (\rho v)_h^- \\
            (\partial_{x_i}\rho_h^-)e_h^- + \rho_h^-(\partial_{x_i} e_h^+) + \partial_{x_i}\left(\rho_h^-\left[(u_h^-)^2 + (v_h^-)^2\right]\right) \\
        \end{pmatrix},\\
        \partial_{x_i}\partial_{x_j}\bQ_h^+ &=\partial_{x_i}\partial_{x_j}\bQ_h^-.
        \end{aligned}
\end{equation}

\subsubsection{Symmetry plane}
$\bQ_h^+$ is computed by determining the mirror image of the velocity with respect to the wall. That is, the normal component of the velocity vector is negated while the tangential component is kept the same. Also, all other variables are extrapolated from the interior state. Furthermore, $\partial_{x_i}\bQ_h^+$ is computed as the velocity vector and $\partial_{x_i}\partial_{x_j}\bQ_h^+$ is equal to the interior second derivative.

\begin{equation}\label{BC:symmetry_plane}
        \bQ^+_h =
        \begin{pmatrix}
            \rho_h^- \\
            \rho_h^-(u_h^--2(\bu_h^-\cdot\bn) n_1) \\
            \rho_h^-(v_h^--2(\bu_h^-\cdot\bn) n_2) \\
             E_h^- \\
        \end{pmatrix},\quad
        \partial_{x_i}\bQ_h^+ =\partial_{x_i}\bQ_h^--2(\partial_n\bQ_h^-\cdot\bn)n_i,\quad
        \partial_{x_i}\partial_{x_j}\bQ_h^+ =\partial_{x_i}\partial_{x_j}\bQ_h^-.
\end{equation}

\section{Numerical Examples}\label{Sec_Num_examples}
In this section, several numerical examples are explored to assess the performance of the proposed DDGIC method. These examples demonstrate the capability of the new DDGIC scheme for achieving optimal order of accuracy as well as capturing correct physical results. In all examples, strong stability preserving (SSP) explicit third order Runge-Kutta scheme in \cite{Shu-Osher-1988,Gottlieb-Shu-Tabmor-2001} is used to evolve the solution polynomial coefficients (\ref{Eqn_DDGIC_ns2d}) in time with the following CFL (Courant-Friedrichs-Levy) condition
\begin{equation}
\Delta t\max_K\left\{\max\left\{\frac{a_K+\lVert\textbf{u}_K\lVert}{h},\frac{\mu_K}{h^2}\right\}\right\}<\omega\lambda,
\end{equation}
where $\lambda$ is the CFL number, $h=\min_K{h_K}$ is the minimum cell diameter in the computational domain, and $\omega$ is the minimum quadrature weight of volume integration. In addition, $a_K$ is the maximum speed of sound, $\lVert\textbf{u}_K\lVert$ is the Euclidean norm of the velocity vector, and $\mu_K$ is the maximum dynamic viscosity computed in the cell $K$. In all numerical experiments, quadrature rules are chosen so that both volume and surface integrals are exact up to polynomials of degree $2k+1$ in each cell.

\begin{figure}[!ht]
	\centering
	\begin{subfigure}[b]{0.22\textwidth}
		\centering
		\includegraphics[scale=0.2]{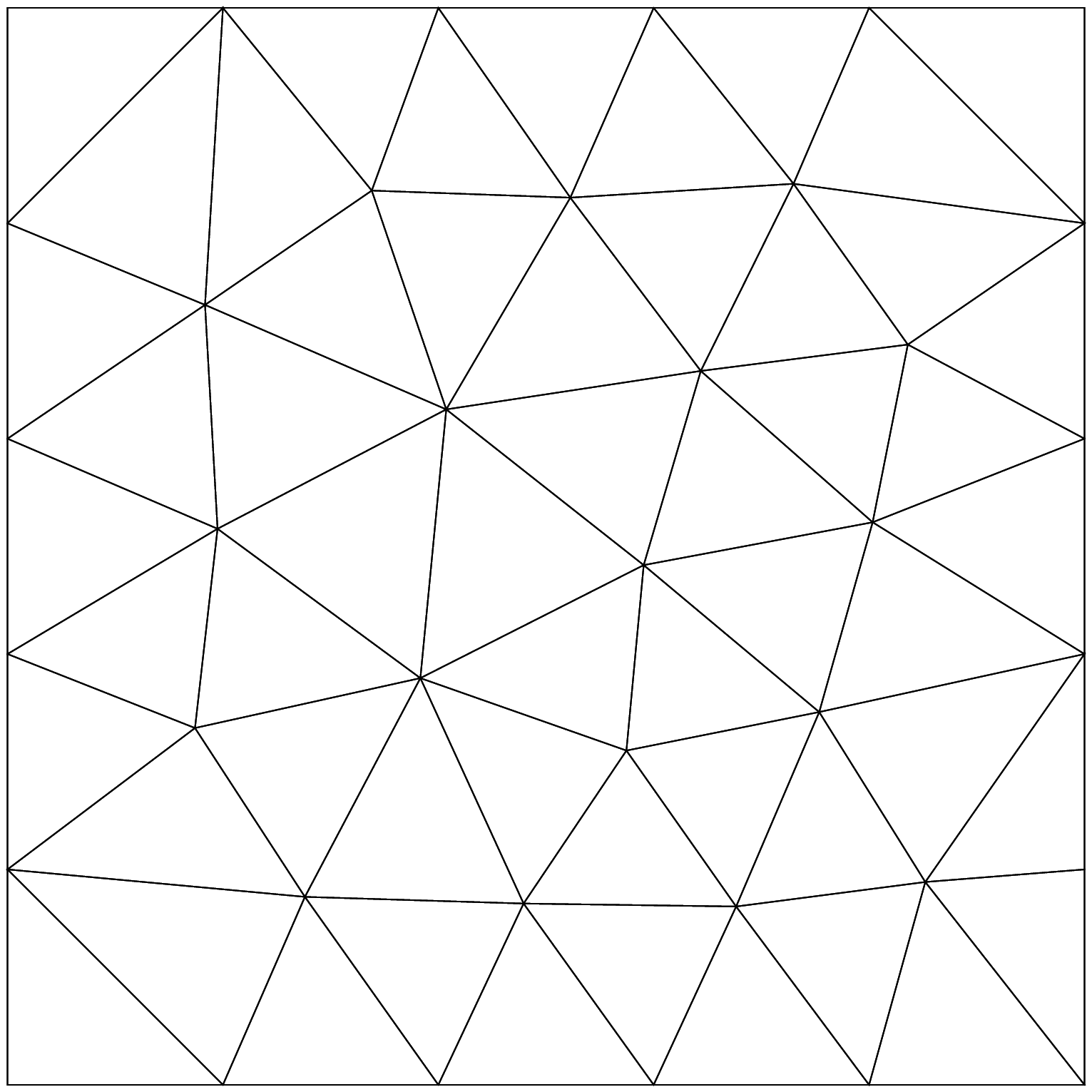}
		\caption{$h$}
	\end{subfigure}
	\begin{subfigure}[b]{0.22\textwidth}
		\centering
		\includegraphics[scale=0.2]{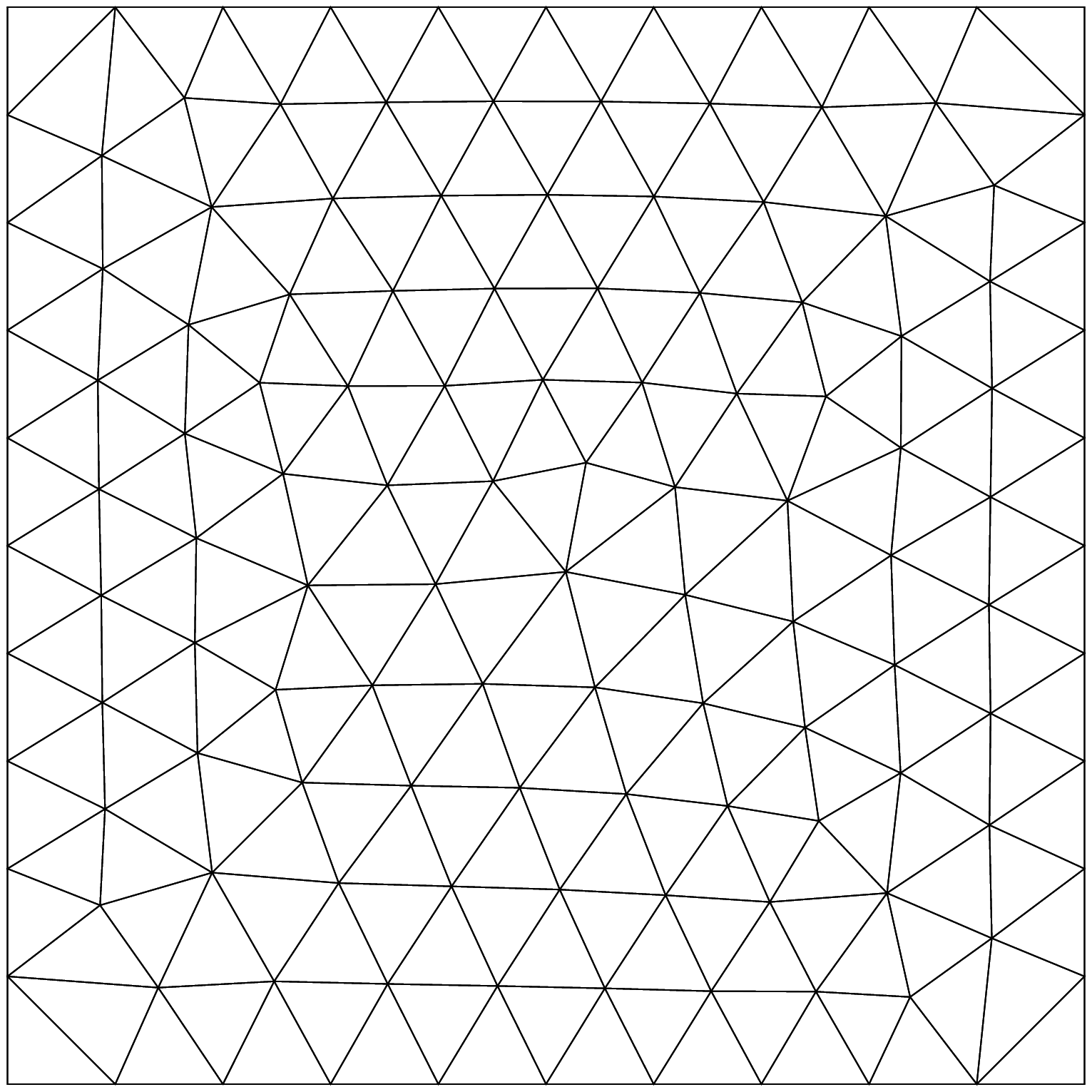}
		\caption{$h/2$}
	\end{subfigure}
	\begin{subfigure}[b]{0.22\textwidth}
		\centering
		\includegraphics[scale=0.2]{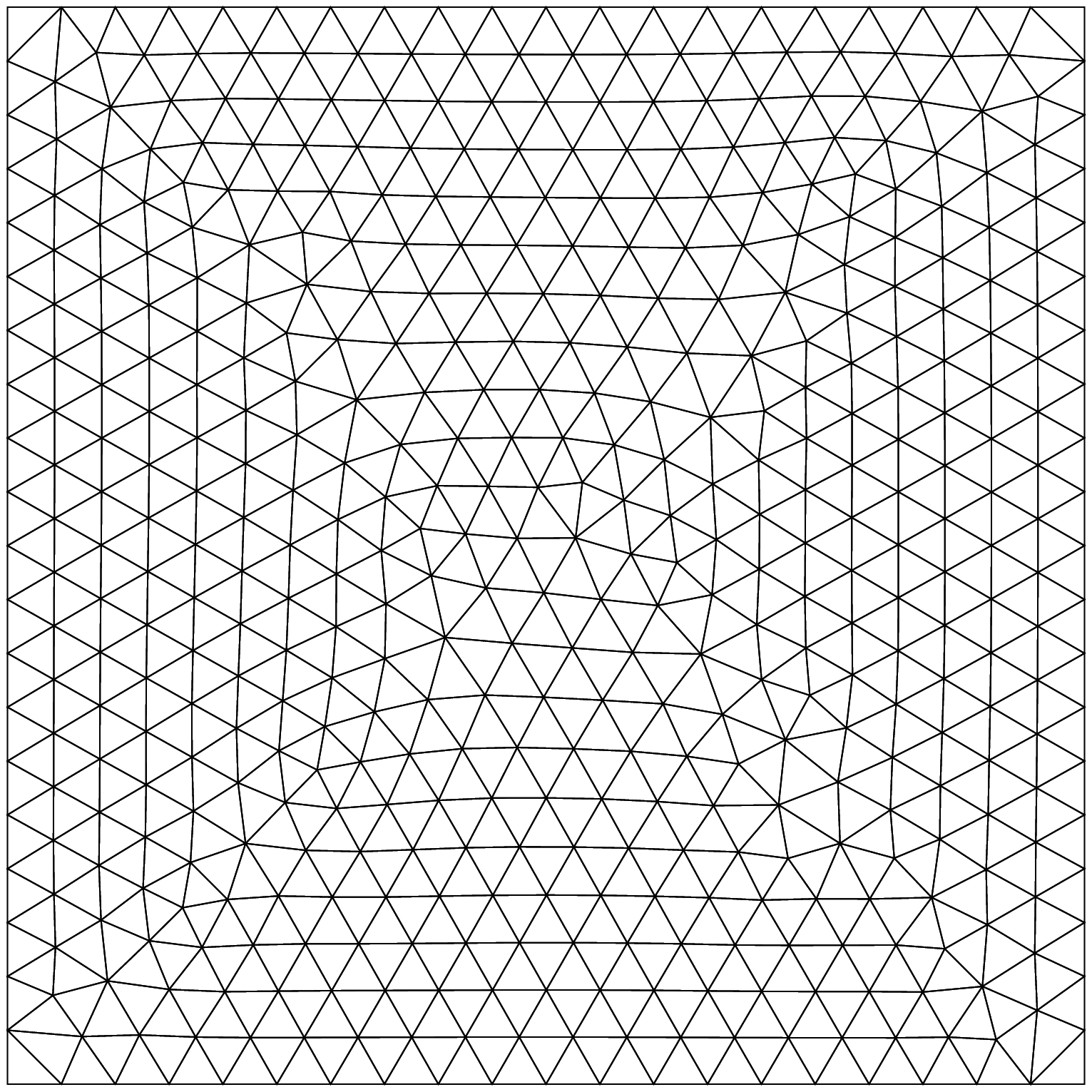}
		\caption{$h/4$}
	\end{subfigure}
	\begin{subfigure}[b]{0.22\textwidth}
		\centering
		\includegraphics[scale=0.2]{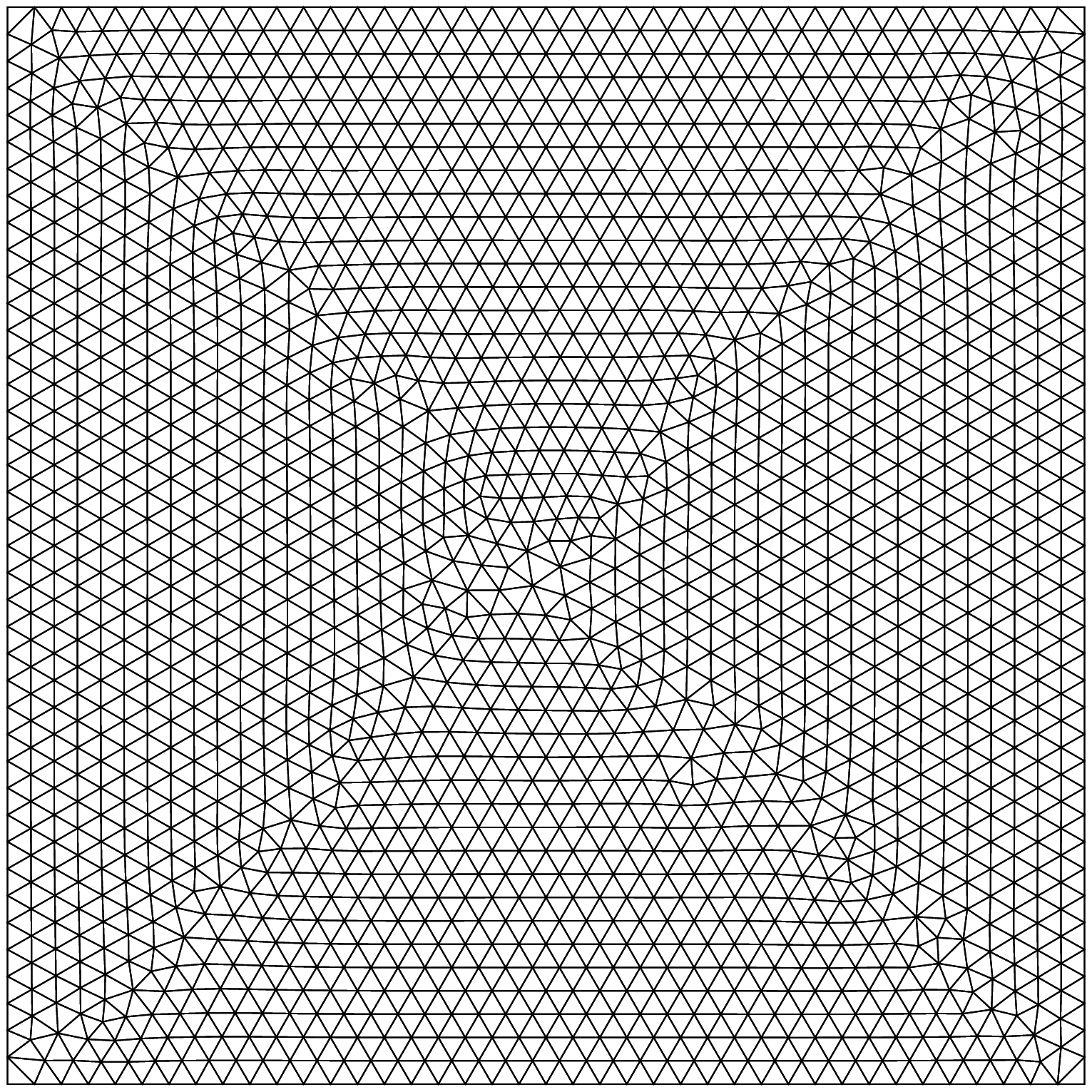}
		\caption{$h/8$}
	\end{subfigure}
	\caption{Uniform unstructured triangular meshes used for all accuracy tests. Here, $h=1/5$ is the length of a typical segment on the domain boundary}\label{Fig_IrregMesh}
\end{figure}

\subsection{Manufactured Solution-1}\label{SubSec_Manuf1}
In this example, we employ the following simple manufactured solution to measure the accuracy of the new DDGIC scheme
\begin{equation*}
\begin{aligned}
\rho(x,y,t)&=1+0.1 \sin(2\pi(x+y)-2t), \\ 
u(x,y,t)&=1+0.1 \sin(2\pi(x+y)-2t), \\
v(x,y,t)&=1+0.1 \cos(2\pi(x+y)-2t), \\
e(x,y,t)&=1+0.1 \cos(2\pi(x+y)-2t). \\
\end{aligned}
\end{equation*} 
The computational domain is set as $\Omega=[0,1]\times[0,1]$, and the set of unstructured triangular meshes used is shown in Figure \ref{Fig_IrregMesh}. Note that these are traveling waves moving in the same direction with the same speed. We use a constant dynamic viscosity of $\mu=10^{-3}$. The CFL number is taken as $\lambda=0.1$ and we run the simulation to final time $T=2\pi$. Since the manufactured solution is periodic in the x- and y- directions, respectively, the periodic boundary conditions \eqref{BC:periodic} are enforced in these directions.

We compute the $L_2$ and $L_\infty$ errors for all conserved variables in this example. The $L_2$ error is computed by a quadrature rule of order $k+1$ while the $L_\infty$ error is computed using 361 points in each cell. This error calculation formula is also employed in the following two numerical examples. To save space, we only list the $L_\infty$ errors and orders in Table \ref{Table_manuf1_Linf}. In this example, an optimal $(k+1)^{th}$ order of accuracy is obtained for each conserved variable with $\mathbb{P}_k$ ($k=1,\cdots,4$) polynomials.

\begin{table}[h!]
\begin{center}
\caption{$L_\infty$ errors for Example \ref{SubSec_Manuf1} at $T=2\pi$}\label{Table_manuf1_Linf}
\begin{tabular}{cccccccc}
\hlineB{3}
\multicolumn{1}{l}{} & \multicolumn{7}{l}{$L_\infty$ errors and orders for $\rho_h$}               \\ \cline{2-8} 
                     & $h$     & $h/2$   & Order & $h/4$   & Order & $h/8$   & Order \\ \hline
$k=1$                & 6.60E-2 & 1.73E-2 & 1.93  & 4.97E-3 & 1.80  & 1.21E-3 & 2.04  \\
$k=2$                & 2.01E-2 & 3.77E-3 & 2.42  & 4.83E-4 & 2.96  & 6.40E-5 & 2.92  \\
$k=3$                & 2.49E-3 & 1.59E-4 & 3.97  & 1.52E-5 & 3.39  & 7.49E-7 & 4.34  \\
$k=4$                & 5.01E-4 & 3.10E-5 & 4.01  & 1.09E-6 & 4.83  & 4.03E-8 & 4.76  \\ \hlineB{2}
\multicolumn{8}{l}{}                                                                 \\ \hlineB{2}
\multicolumn{1}{l}{} & \multicolumn{7}{l}{$L_\infty$ errors and orders for $(\rho u)_h$}             \\ \cline{2-8} 
                     & $h$     & $h/2$   & Order & $h/4$   & Order & $h/8$   & Order \\ \hline
$k=1$                & 1.33E-1 & 3.23E-2 & 2.04  & 1.01E-2 & 1.68  & 2.24E-3 & 2.17  \\
$k=2$                & 3.66E-2 & 7.23E-3 & 2.34  & 9.05E-4 & 3.00  & 1.13E-4 & 3.01  \\
$k=3$                & 6.03E-3 & 3.45E-4 & 4.13  & 3.31E-5 & 3.38  & 1.40E-6 & 4.56  \\
$k=4$                & 8.48E-4 & 5.32E-5 & 3.99  & 1.74E-6 & 4.93  & 5.66E-8 & 4.94  \\ \hlineB{2}
\multicolumn{8}{l}{}                                                                 \\ \hlineB{2}
\multicolumn{1}{l}{} & \multicolumn{7}{l}{$L_\infty$ errors and orders for $(\rho v)_h$}             \\ \cline{2-8} 
                     & $h$     & $h/2$   & Order & $h/4$   & Order & $h/8$   & Order \\ \hline
$k=1$                & 1.05E-1 & 3.20E-2 & 1.72  & 7.68E-3 & 2.06  & 1.78E-3 & 2.11  \\
$k=2$                & 2.94E-2 & 5.59E-3 & 2.39  & 7.10E-4 & 2.98  & 9.01E-5 & 2.98  \\
$k=3$                & 4.73E-3 & 4.18E-4 & 3.50  & 2.11E-5 & 4.31  & 1.22E-6 & 4.11  \\
$k=4$                & 1.23E-3 & 6.99E-5 & 4.14  & 2.23E-6 & 4.97  & 7.03E-8 & 4.99  \\ \hlineB{2}
\multicolumn{8}{l}{}                                                                 \\ \hlineB{2}
\multicolumn{1}{l}{} & \multicolumn{7}{l}{$L_\infty$ errors and orders for $E_h$}                  \\ \cline{2-8} 
                     & $h$     & $h/2$   & Order & $h/4$   & Order & $h/8$   & Order \\ \hline
$k=1$                & 2.75E-1 & 7.37E-2 & 1.90  & 1.70E-2 & 2.12  & 4.04E-3 & 2.07  \\
$k=2$                & 7.37E-2 & 1.40E-2 & 2.39  & 1.78E-3 & 2.98  & 2.25E-4 & 2.99  \\
$k=3$                & 1.32E-2 & 8.82E-4 & 3.90  & 4.23E-5 & 4.38  & 2.48E-6 & 4.09  \\
$k=4$                & 2.72E-3 & 1.57E-4 & 4.12  & 4.97E-6 & 4.98  & 1.54E-7 & 5.02  \\ \hlineB{3}
\end{tabular}
\end{center}
\end{table}

\subsection{Manufactured Solution-2}\label{SubSec_Manuf2}
In this example, the new DDGIC method is tested with a more complicated manufactured solution. We now consider a wave-packet solution:
\begin{equation*}
\begin{aligned}
\rho(x,y,t)&=1- 0.1 \sin(4\pi x + 4\pi t) \cos(2\pi y - 2\pi t), \\ 
u(x,y,t)&=2   + 0.2 \sin(2\pi x - 2\pi t) \cos(4\pi y - 4\pi t), \\
v(x,y,t)&=3   + 0.3 \cos(2\pi x - 2\pi t) \sin(4\pi y + 4\pi t), \\
e(x,y,t)&=50  - 10  \cos(2\pi x - 4\pi t) \sin(4\pi y + 4\pi t). 
\end{aligned}
\end{equation*} 
Note that these waves are traveling in different directions with different phase and group speeds. The CFL number is taken as $\lambda=0.1$, the dynamic viscosity is $\mu=10^{-2}$, and the final time is $T=1$. The periodic boundary conditions \eqref{BC:periodic} are employed in x- and y- directions. The same set of meshes is used as in the previous example. In Table \ref{Table_manuf2_L2}, we only list the $L_2$ errors and orders to save space. Again, we obtain an optimal $(k+1)^{th}$ order of accuracy in both $L_2$ and $L_\infty$ sense for all conserved variables. 

\begin{table}[h!]
\begin{center}
\caption{$L_2$ Errors for Example \ref{SubSec_Manuf2} at $T=1$}\label{Table_manuf2_L2}
\begin{tabular}{cccccccc}
\hlineB{3}
\multicolumn{1}{l}{} & \multicolumn{7}{l}{$L_2$ errors and orders for $\rho_h$}                    \\ \cline{2-8} 
                     & $h$     & $h/2$   & Order & $h/4$   & Order & $h/8$   & Order \\ \hline
$k=1$                & 4.16E-2 & 5.34E-3 & 2.96  & 1.01E-3 & 2.40  & 2.22E-4 & 2.19  \\
$k=2$                & 5.81E-3 & 9.96E-4 & 2.54  & 1.49E-4 & 2.75  & 1.59E-5 & 3.22  \\
$k=3$                & 1.63E-3 & 8.75E-5 & 4.22  & 4.57E-6 & 4.26  & 2.54E-7 & 4.17  \\
$k=4$                & 2.55E-4 & 8.83E-6 & 4.85  & 2.43E-7 & 5.18  & 6.58E-9 & 5.21  \\ \hlineB{2}
\multicolumn{8}{l}{}                                                                 \\ \hlineB{2}
\multicolumn{1}{l}{} & \multicolumn{7}{l}{$L_2$ errors and orders for $(\rho u)_h$}                  \\ \cline{2-8} 
                     & $h$     & $h/2$   & Order & $h/4$   & Order & $h/8$   & Order \\ \hline
$k=1$                & 1.05E-1 & 1.66E-2 & 2.66  & 3.28E-3 & 2.33  & 7.06E-4 & 2.22  \\
$k=2$                & 1.54E-2 & 2.29E-3 & 2.75  & 3.17E-4 & 2.85  & 3.59E-5 & 3.14  \\
$k=3$                & 3.84E-3 & 1.87E-4 & 4.36  & 9.96E-6 & 4.23  & 5.59E-7 & 4.16  \\
$k=4$                & 5.65E-4 & 1.67E-5 & 5.08  & 5.08E-7 & 5.04  & 1.46E-8 & 5.12  \\ \hlineB{2}
\multicolumn{8}{l}{}                                                                 \\ \hlineB{2}
\multicolumn{1}{l}{} & \multicolumn{7}{l}{$L_2$ errors and orders for $(\rho v)_h$}                  \\ \cline{2-8} 
                     & $h$     & $h/2$   & Order & $h/4$   & Order & $h/8$   & Order \\ \hline
$k=1$                & 9.12E-2 & 1.71E-2 & 2.41  & 3.53E-3 & 2.28  & 7.61E-4 & 2.22  \\
$k=2$                & 1.61E-2 & 2.07E-3 & 2.96  & 2.45E-4 & 3.08  & 2.87E-5 & 3.10  \\
$k=3$                & 3.01E-3 & 1.67E-4 & 4.17  & 9.47E-6 & 4.14  & 5.37E-7 & 4.14  \\
$k=4$                & 5.64E-4 & 1.49E-5 & 5.24  & 4.66E-7 & 5.00  & 1.51E-8 & 4.95  \\ \hlineB{2}
\multicolumn{8}{l}{}                                                                 \\ \hlineB{2}
\multicolumn{1}{l}{} & \multicolumn{7}{l}{$L_2$ errors and orders for $E_h$}                       \\ \cline{2-8} 
                     & $h$     & $h/2$   & Order & $h/4$   & Order & $h/8$   & Order \\ \hline
$k=1$                & 1.61E+0 & 3.06E-1 & 2.39  & 7.03E-2 & 2.12  & 1.70E-2 & 2.05  \\
$k=2$                & 2.96E-1 & 4.36E-2 & 2.76  & 6.39E-3 & 2.77  & 8.62E-4 & 2.89  \\
$k=3$                & 6.59E-2 & 3.76E-3 & 4.13  & 2.20E-4 & 4.09  & 1.41E-5 & 3.97  \\
$k=4$                & 1.08E-2 & 3.83E-4 & 4.81  & 1.24E-5 & 4.95  & 4.01E-7 & 4.95  \\ \hlineB{3}
\end{tabular}
\end{center}
\end{table}

\subsection{A pressure pulse in a periodic domain}\label{SubSec_PressurePulse}
In this example, we consider the Test Case II of \cite{Hesthaven_nodal}. The initial conditions are given as
\begin{equation}
\begin{aligned}
\rho(x,y,0) &= 1, \\
u(x,y,0)    &= v(x,y,0) = 0, \\
E(x,y,0)    &= \frac{12}{\gamma-1} + \frac{1}{2}\exp{\left(-\left( \cos(\pi x)^2 + \cos(\pi y)^2\right)\right)}.
\end{aligned}
\end{equation}
The reference solution is generated with high order $5^{th}$ degree  polynomial approximation and on a refined mesh of $h/16$. The CFL number is taken as  $\lambda=0.1$, the dynamic viscosity as $\mu=10^{-2}$, and the final time as $T=0.1$. Periodic boundary conditions \eqref{BC:periodic} are enforced in x- and y- directions. The set of unstructured triangular meshes used in this example is the same as in Example \ref{SubSec_Manuf1}. The $L_2$ and $L_\infty$ errors and orders are given in Table \ref{Table_pressure_pulse_L2} and \ref{Table_pressure_pulse_Linf}. As in the previous cases, the new DDGIC method achieves the expected $(k+1)^{th}$ optimal order of accuracy. 

\begin{table}[h!]
\begin{center}
\caption{$L_2$ Errors for Example \ref{SubSec_PressurePulse} at $T=0.1$}\label{Table_pressure_pulse_L2}
\begin{tabular}{cccccccc}
\hlineB{3}
\multicolumn{1}{l}{} & \multicolumn{7}{l}{$L_2$ errors and orders for $\rho_h$}                     \\ \cline{2-8} 
                     & $h$     & $h/2$   & Order & $h/4$   & Order & $h/8$    & Order \\ \hline
$k=1$                & 6.27E-4 & 1.28E-4 & 2.30  & 2.64E-5 & 2.27  & 5.60E-6  & 2.24  \\
$k=2$                & 1.01E-4 & 1.46E-5 & 2.79  & 1.91E-6 & 2.93  & 2.19E-7  & 3.12  \\
$k=3$                & 3.21E-5 & 1.26E-6 & 4.67  & 7.75E-8 & 4.03  & 4.27E-9  & 4.18  \\
$k=4$                & 4.79E-6 & 1.59E-7 & 4.92  & 4.49E-9 & 5.14  & 1.50E-10 & 4.90  \\ \hlineB{2}
\multicolumn{8}{l}{}                                                                  \\ \hlineB{2}
\multicolumn{1}{l}{} & \multicolumn{7}{l}{$L_2$ errors and orders for $(\rho u)_h$}                   \\ \cline{2-8} 
                     & $h$     & $h/2$   & Order & $h/4$   & Order & $h/8$    & Order \\ \hline
$k=1$                & 1.50E-3 & 3.41E-4 & 2.13  & 6.47E-5 & 2.40  & 1.24E-5  & 2.38  \\
$k=2$                & 3.12E-4 & 2.75E-5 & 3.50  & 3.07E-6 & 3.16  & 3.77E-7  & 3.02  \\
$k=3$                & 6.08E-5 & 2.57E-6 & 4.56  & 1.60E-7 & 4.00  & 9.28E-9  & 4.11  \\
$k=4$                & 1.05E-5 & 2.62E-7 & 5.33  & 8.53E-9 & 4.94  & 2.76E-10 & 4.95  \\ \hlineB{2}
\multicolumn{8}{l}{}                                                                  \\ \hlineB{2}
\multicolumn{1}{l}{} & \multicolumn{7}{l}{$L_2$ errors and orders for $(\rho v)_h$}                   \\ \cline{2-8} 
                     & $h$     & $h/2$   & Order & $h/4$   & Order & $h/8$    & Order \\ \hline
$k=1$                & 1.35E-3 & 3.49E-4 & 1.95  & 6.65E-5 & 2.39  & 1.25E-5  & 2.41  \\
$k=2$                & 2.62E-4 & 2.77E-5 & 3.24  & 3.11E-6 & 3.16  & 3.80E-7  & 3.03  \\
$k=3$                & 5.12E-5 & 2.59E-6 & 4.30  & 1.61E-7 & 4.01  & 9.31E-9  & 4.11  \\
$k=4$                & 8.26E-6 & 2.69E-7 & 4.94  & 8.65E-9 & 4.96  & 2.76E-10 & 4.97  \\ \hlineB{2}
\multicolumn{8}{l}{}                                                                  \\ \hlineB{2}
\multicolumn{1}{l}{} & \multicolumn{7}{l}{$L_2$ errors and orders for $E_h$}                        \\ \cline{2-8} 
                     & $h$     & $h/2$   & Order & $h/4$   & Order & $h/8$    & Order \\ \hline
$k=1$                & 1.51E-2 & 3.61E-3 & 2.06  & 7.66E-4 & 2.24  & 1.56E-4  & 2.30  \\
$k=2$                & 3.28E-3 & 3.31E-4 & 3.31  & 3.23E-5 & 3.36  & 3.68E-6  & 3.13  \\
$k=3$                & 7.71E-4 & 2.66E-5 & 4.86  & 1.65E-6 & 4.02  & 9.88E-8  & 4.06  \\
$k=4$                & 1.16E-4 & 2.91E-6 & 5.31  & 7.91E-8 & 5.20  & 3.12E-9  & 4.66  \\ \hlineB{3}
\end{tabular}
\end{center}
\end{table}
\newpage

\begin{table}[h!]
\begin{center}
\caption{$L_\infty$ Errors for Example \ref{SubSec_PressurePulse} at $T=0.1$}\label{Table_pressure_pulse_Linf}
\begin{tabular}{cccccccc}
\hlineB{3}
\multicolumn{1}{l}{} & \multicolumn{7}{l}{$L_\infty$ errors and orders for $\rho_h$}               \\ \cline{2-8} 
                     & $h$     & $h/2$   & Order & $h/4$   & Order & $h/8$   & Order \\ \hline
$k=1$                & 3.10E-3 & 1.36E-3 & 1.19  & 3.51E-4 & 1.96  & 9.86E-5 & 1.83  \\
$k=2$                & 6.86E-4 & 2.49E-4 & 1.46  & 2.24E-5 & 3.47  & 2.78E-6 & 3.01  \\
$k=3$                & 3.24E-4 & 2.18E-5 & 3.89  & 2.26E-6 & 3.27  & 2.02E-7 & 3.48  \\
$k=4$                & 5.98E-5 & 5.79E-6 & 3.37  & 1.18E-7 & 5.62  & 3.70E-9 & 4.99  \\ \hlineB{2}
\multicolumn{8}{l}{}                                                                 \\ \hlineB{2}
\multicolumn{1}{l}{} & \multicolumn{7}{l}{$L_\infty$ errors and orders for $(\rho u)_h$}             \\ \cline{2-8} 
                     & $h$     & $h/2$   & Order & $h/4$   & Order & $h/8$   & Order \\ \hline
$k=1$                & 4.59E-3 & 1.17E-3 & 1.97  & 3.49E-4 & 1.75  & 1.16E-4 & 1.59  \\
$k=2$                & 1.25E-3 & 2.98E-4 & 2.07  & 4.33E-5 & 2.78  & 4.91E-6 & 3.14  \\
$k=3$                & 6.86E-4 & 6.38E-5 & 3.43  & 4.15E-6 & 3.94  & 2.55E-7 & 4.02  \\
$k=4$                & 1.49E-4 & 6.63E-6 & 4.49  & 1.81E-7 & 5.19  & 5.18E-9 & 5.13  \\ \hlineB{2}
\multicolumn{8}{l}{}                                                                 \\ \hlineB{2}
\multicolumn{1}{l}{} & \multicolumn{7}{l}{$L_\infty$ errors and orders for $(\rho v)_h$}             \\ \cline{2-8} 
                     & $h$     & $h/2$   & Order & $h/4$   & Order & $h/8$   & Order \\ \hline
$k=1$                & 3.94E-3 & 1.20E-3 & 1.71  & 3.33E-4 & 1.85  & 1.18E-4 & 1.49  \\
$k=2$                & 1.12E-3 & 2.84E-4 & 1.98  & 4.70E-5 & 2.59  & 5.28E-6 & 3.15  \\
$k=3$                & 6.52E-4 & 5.23E-5 & 3.64  & 3.40E-6 & 3.95  & 2.81E-7 & 3.60  \\
$k=4$                & 1.06E-4 & 6.56E-6 & 4.02  & 2.30E-7 & 4.83  & 5.16E-9 & 5.48  \\ \hlineB{2}
\multicolumn{8}{l}{}                                                                 \\ \hlineB{2}
\multicolumn{1}{l}{} & \multicolumn{7}{l}{$L_\infty$ errors and orders for $E_h$}                  \\ \cline{2-8} 
                     & $h$     & $h/2$   & Order & $h/4$   & Order & $h/8$   & Order \\ \hline
$k=1$                & 4.88E-2 & 1.24E-2 & 1.98  & 3.99E-3 & 1.64  & 1.44E-3 & 1.47  \\
$k=2$                & 1.40E-2 & 3.06E-3 & 2.20  & 6.27E-4 & 2.29  & 5.99E-5 & 3.39  \\
$k=3$                & 7.24E-3 & 6.02E-4 & 3.59  & 4.44E-5 & 3.76  & 3.39E-6 & 3.71  \\
$k=4$                & 1.62E-3 & 6.72E-5 & 4.59  & 3.29E-6 & 4.35  & 7.39E-8 & 5.48  \\ \hlineB{3}
\end{tabular}
\end{center}
\end{table}

\subsection{Steady flow over a flat plate}\label{SubSec_FlatPLate}
\begin{figure}[!ht]
\begin{center}
	\includegraphics[scale=0.6]{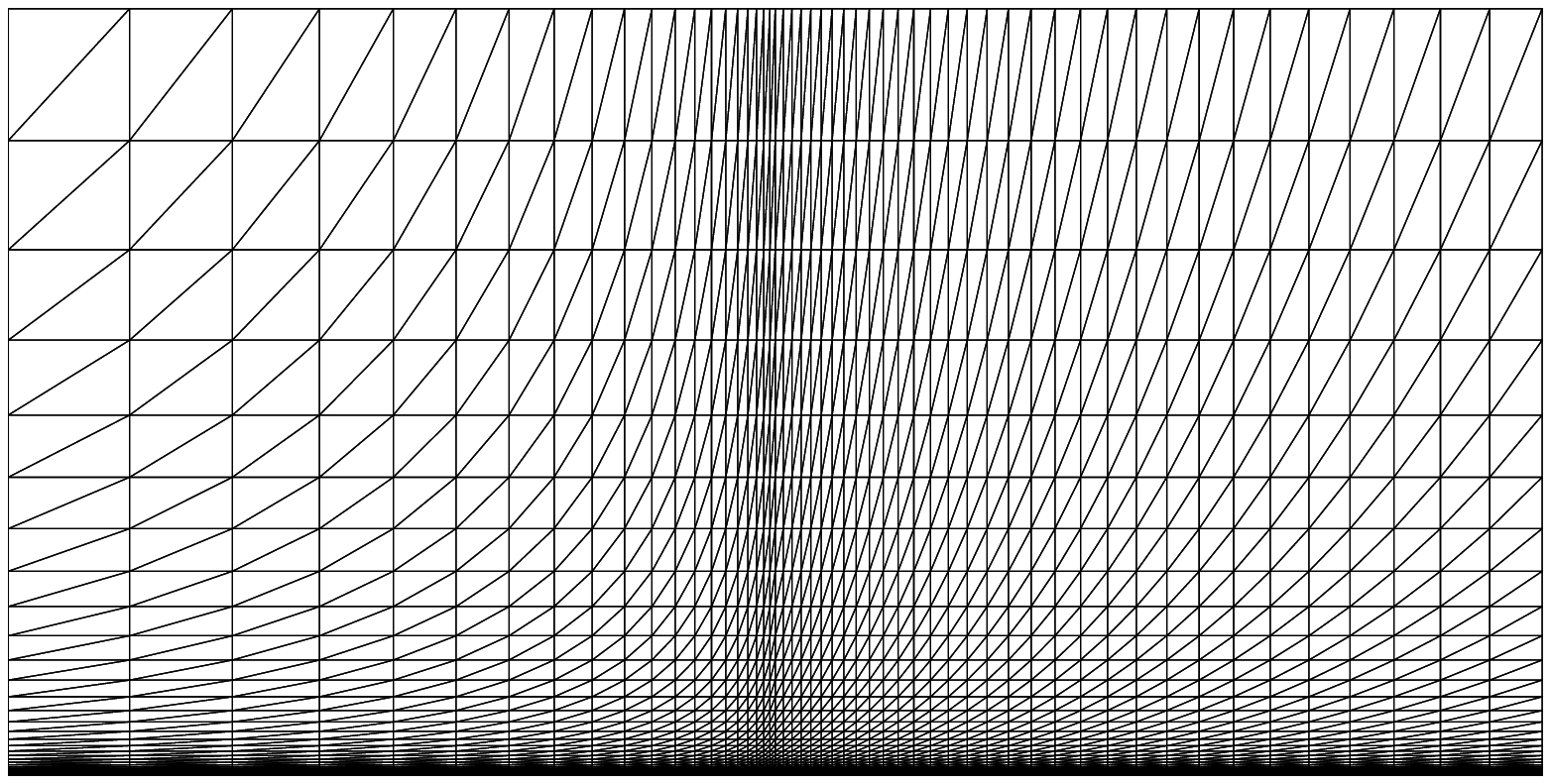}
	\caption{The unstructured triangular mesh used in Example (\ref{SubSec_FlatPLate})}\label{Fig_Blasius_mesh}
\end{center}
\end{figure}

In this example, the simulation setup is similar to example 4.2 of \cite{KANNAN20102007}. We consider a laminar compressible flow over a flat plate with a Reynolds number of $Re=10^4$ based on the inflow reference conditions and the plate length. The inflow Mach number is $M_{\infty}=0.3$. The computational domain is given as $\Omega=[-1,1]\times[0,1]$. We apply the inflow and farfield boundary condition \eqref{BC:inflow_farfield} at the inlet ($x=-1$) and the top  ($y=1$) planes, respectively. Moreover, the outflow boundary condition \eqref{BC:outflow} is enforced at the exit plane ($x=1$). At $y=0$, the symmetry-plane boundary condition \eqref{BC:symmetry_plane} is applied for $-1<x<0$ while the adiabatic viscous wall boundary condition \eqref{BC:adiabatic_wall} is employed for $0<x<1$. We denote the latter as the flat plate. The mesh consists of 3368 elements. Note that 32 cells are located on the flat plate and we have 20 cells in the vertical direction, see Figure \ref{Fig_Blasius_mesh}. Also, the flow field is assumed to reach a steady-state after the 2-norm of the residual of each equation is below $10^{-5}$. 

\begin{figure}[!ht]
	\centering
	\begin{subfigure}[b]{0.45\textwidth}
		\centering
		\includegraphics[scale=0.6]{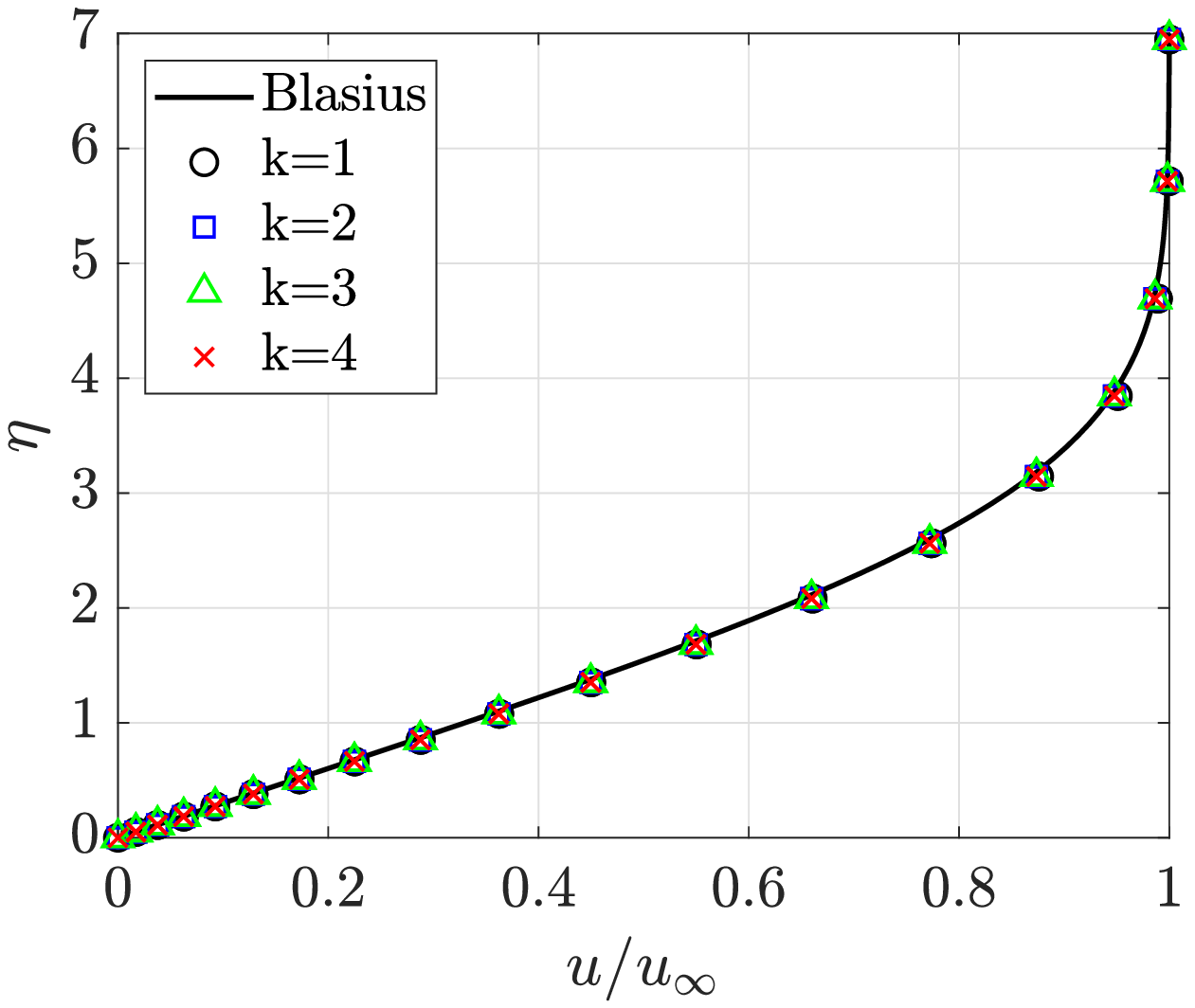}
		\caption{}
	\end{subfigure}
	\begin{subfigure}[b]{0.45\textwidth}
		\centering
		\includegraphics[scale=0.6]{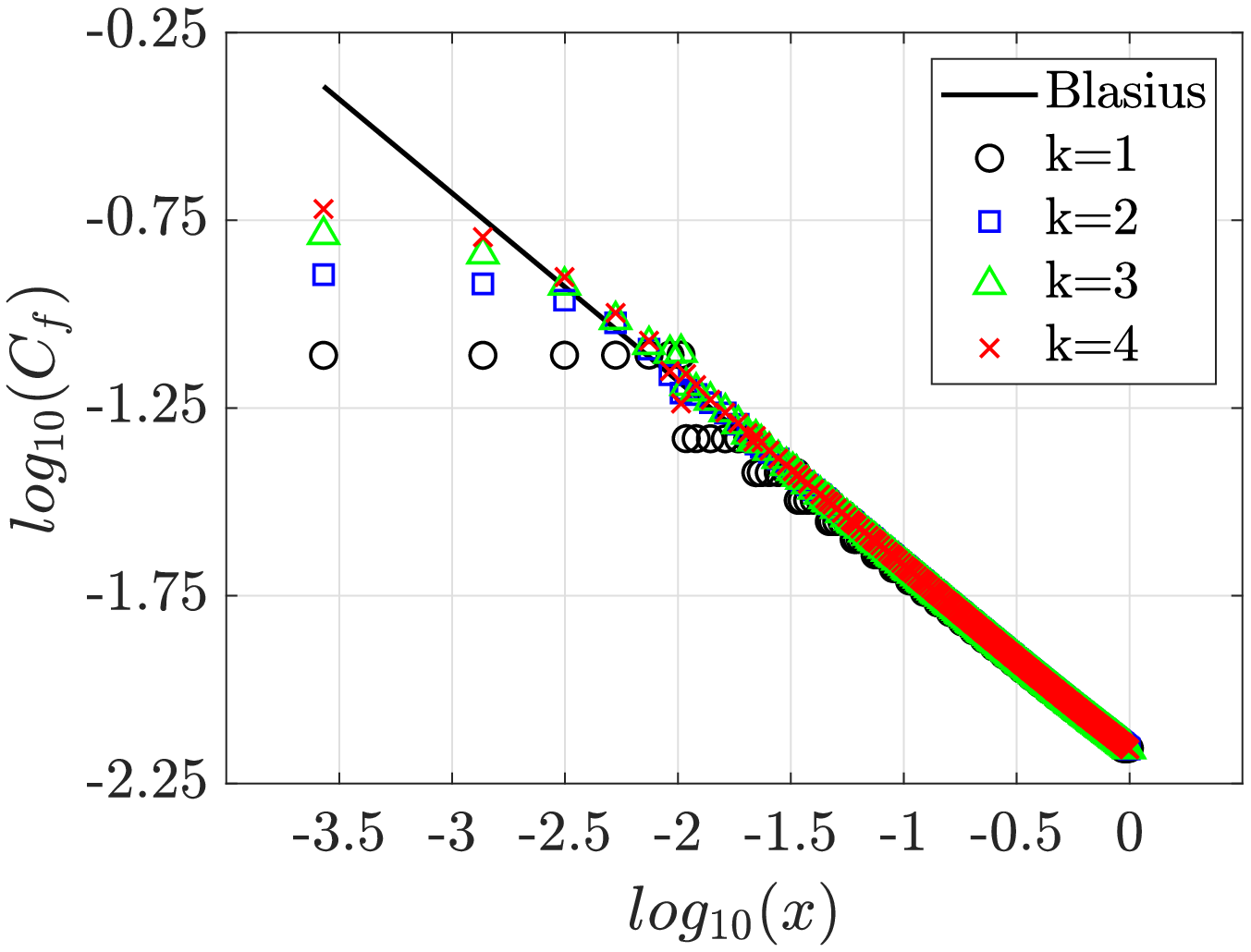}
		\caption{}
	\end{subfigure}
	\caption{The solutions obtained with the new DDGIC method in Example \ref{SubSec_FlatPLate} (a) u-velocity profile at $x=1$, and (b) the skin friction coefficient }\label{Fig_blasius}
\end{figure}

Note that for this inflow Mach number, flow is nearly incompressible so that incompressible laminar boundary layer theory can be used to assess the performance of our method. In Figure \ref{Fig_blasius}a, the DDGIC results of nondimensional u-velocity profile for several polynomial degrees $k$ at the exit plane are compared to the velocity profile of the incompressible Blasius flow. Results corresponding to DDGIC solutions are calculated at the center of each triangle face lying on the exit plane. Note that the wall-normal direction is denoted by the standard nondimensional variable $\eta(x)=y/\delta(x)$ where $\delta(x)=\sqrt{\mu x/\rho_{\infty}u_{\infty}}$. Overall, we observe an excellent agreement in all velocity profiles.

Furthermore, the DDGIC results for skin friction coefficient values along the plane are compared to that obtained from the incompressible theory in Figure \ref{Fig_blasius}b. In this case, DDGIC results are calculated at 20 points for each triangle face lying on the no-slip boundary. The DDGIC results and incompressible theory seem almost identical on a linear plot, thus, the comparison of results is made on a logarithmic plot. It is observed that there is a discrepancy in the skin friction coefficient at the leading edge of the plate. Since the skin friction coefficient is based on the wall-normal derivative of u-velocity, the second order ($k=1$) solution is almost linear in each cell, and thus has the highest discrepancy as seen in Figure \ref{Fig_blasius}b. It is also clear that the discrepancy gets smaller as further p-refinement is applied. This is, in fact, the expected behaviour near the leading edge and further improvement is also possible with h-refinement.

\subsection{Steady flow over a cylinder}\label{SubSec_SteadyCylinder}
\begin{figure}[!ht]
\begin{center}
	\includegraphics[scale=0.5]{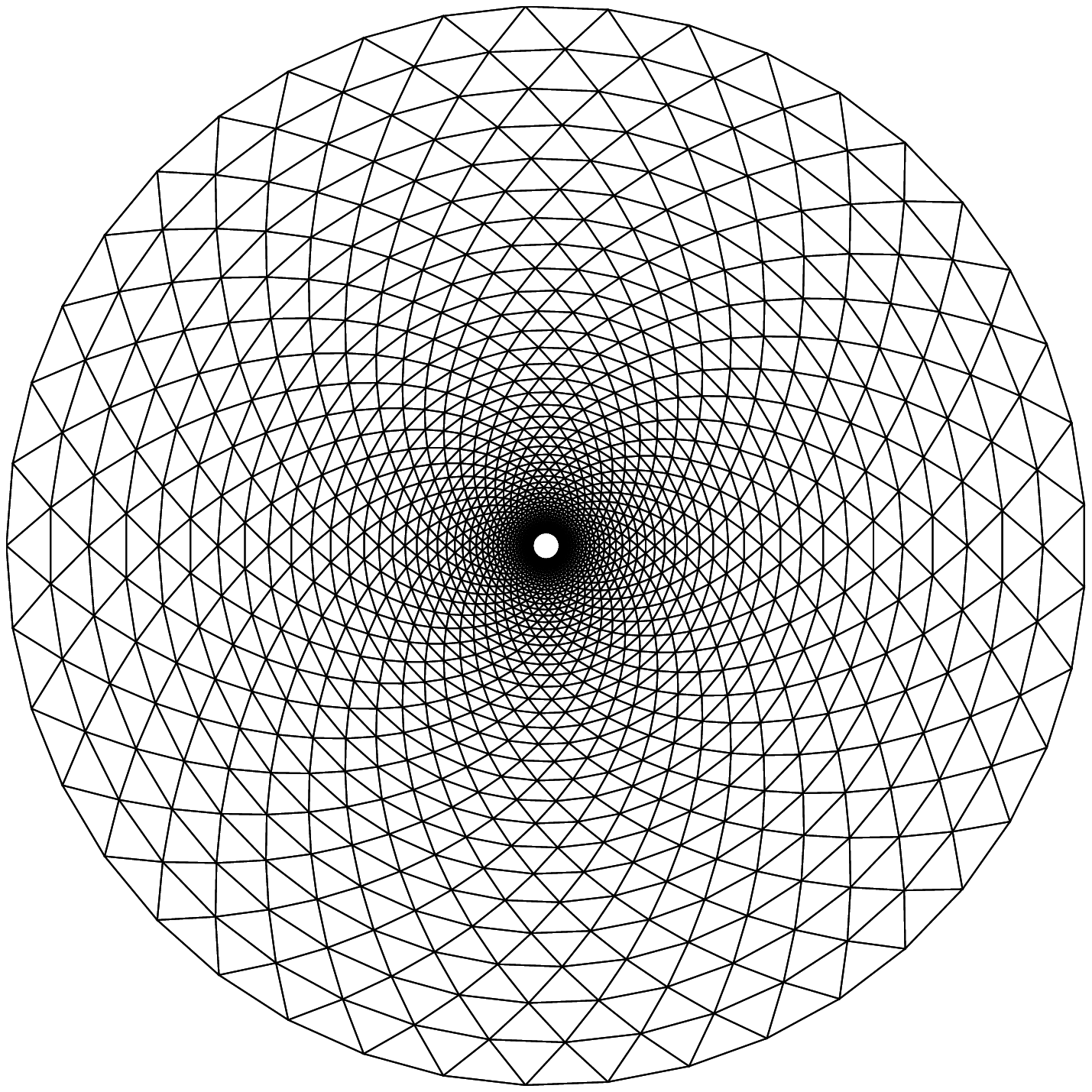}
	\caption{The unstructured triangular mesh used in Example (\ref{SubSec_SteadyCylinder})}\label{Fig_cylinder_mesh}
\end{center}
\end{figure}

In this example, we consider a steady subsonic laminar flow at $Re=40$ based on the inflow reference conditions and the cylinder diameter. The flow enters the computational domain with a Mach number of $M_\infty=0.2$. The computational domain is of a circular shape that extends 20 times the diameter of the cylinder, and consists of $3280$ elements with 41 elements on the cylinder and 40 elements in the radial direction as shown in Figure \ref{Fig_cylinder_mesh}. On the cylinder, we apply the adiabatic viscous wall boundary condition \eqref{BC:adiabatic_wall}. In the remaining parts of the domain boundary, we assume the farfield boundary condition \eqref{BC:inflow_farfield}. Moreover, the flow field is assumed to reach a steady-state after 2-norm of the residual of each equation reaches $10^{-6}$. 

There are several important features of this example used for validation purposes in the literature. It is typical to measure the drag coefficient, which is a nondimensional measure of the drag force is exerted on a body immersed in a fluid. After the flow passes over the cylinder, the boundary layer is said to be \textit{separated} due to viscous and adverse pressure affects. It is thus a common practice to estimate the separation angle and compare it to available reference data. Due to separation, two circulation bubbles occur in the wake of the cylinder. The size of these bubbles in the horizontal and vertical directions as well as the distance between their centers are another important parameters for validating numerical results.

\begin{figure}[!ht]
	\centering
	\begin{subfigure}[b]{0.45\textwidth}
		\centering
		\includegraphics[scale=0.4]{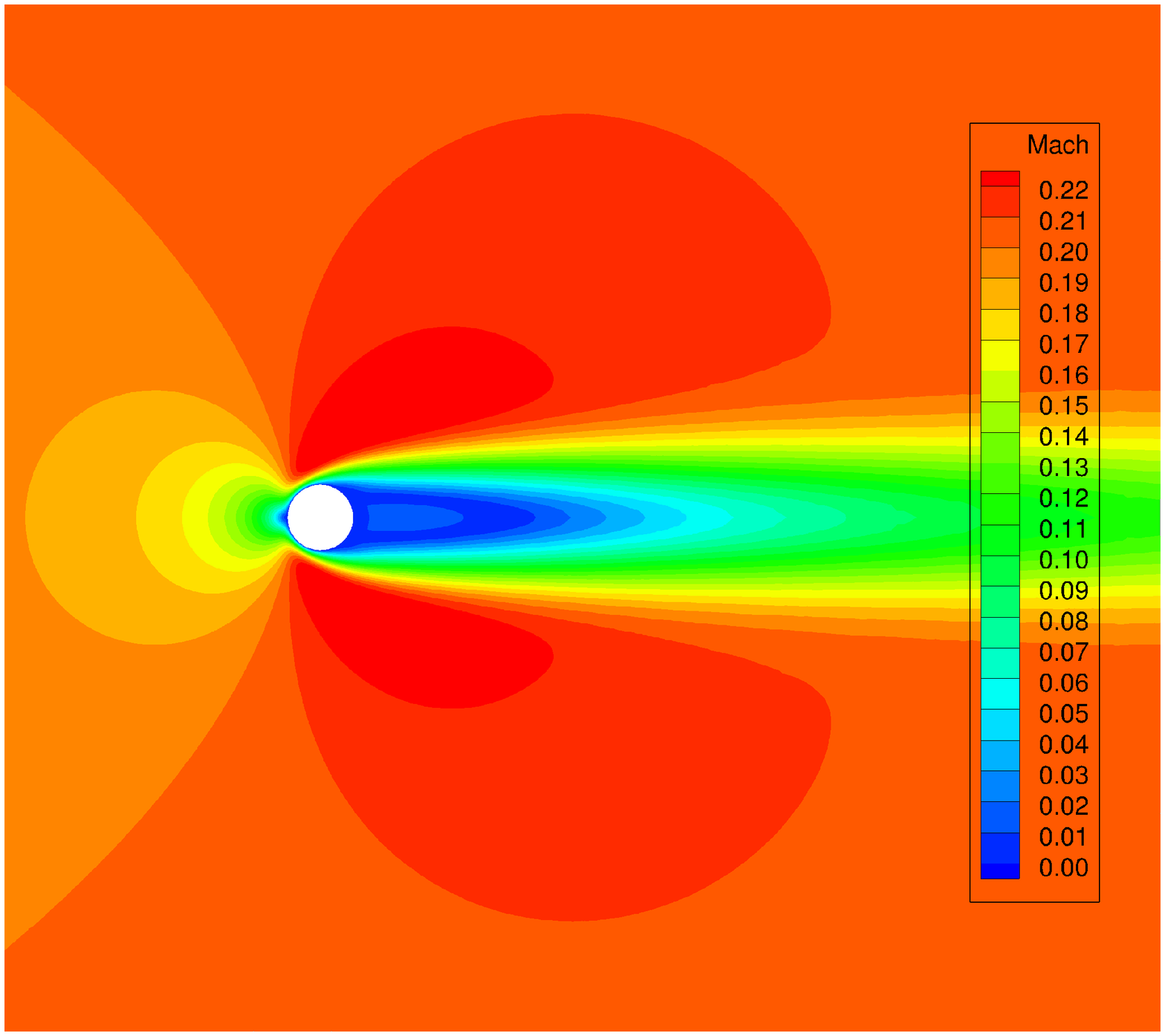}
		\caption{}
	\end{subfigure}
	\begin{subfigure}[b]{0.45\textwidth}
		\centering
		\includegraphics[scale=0.4]{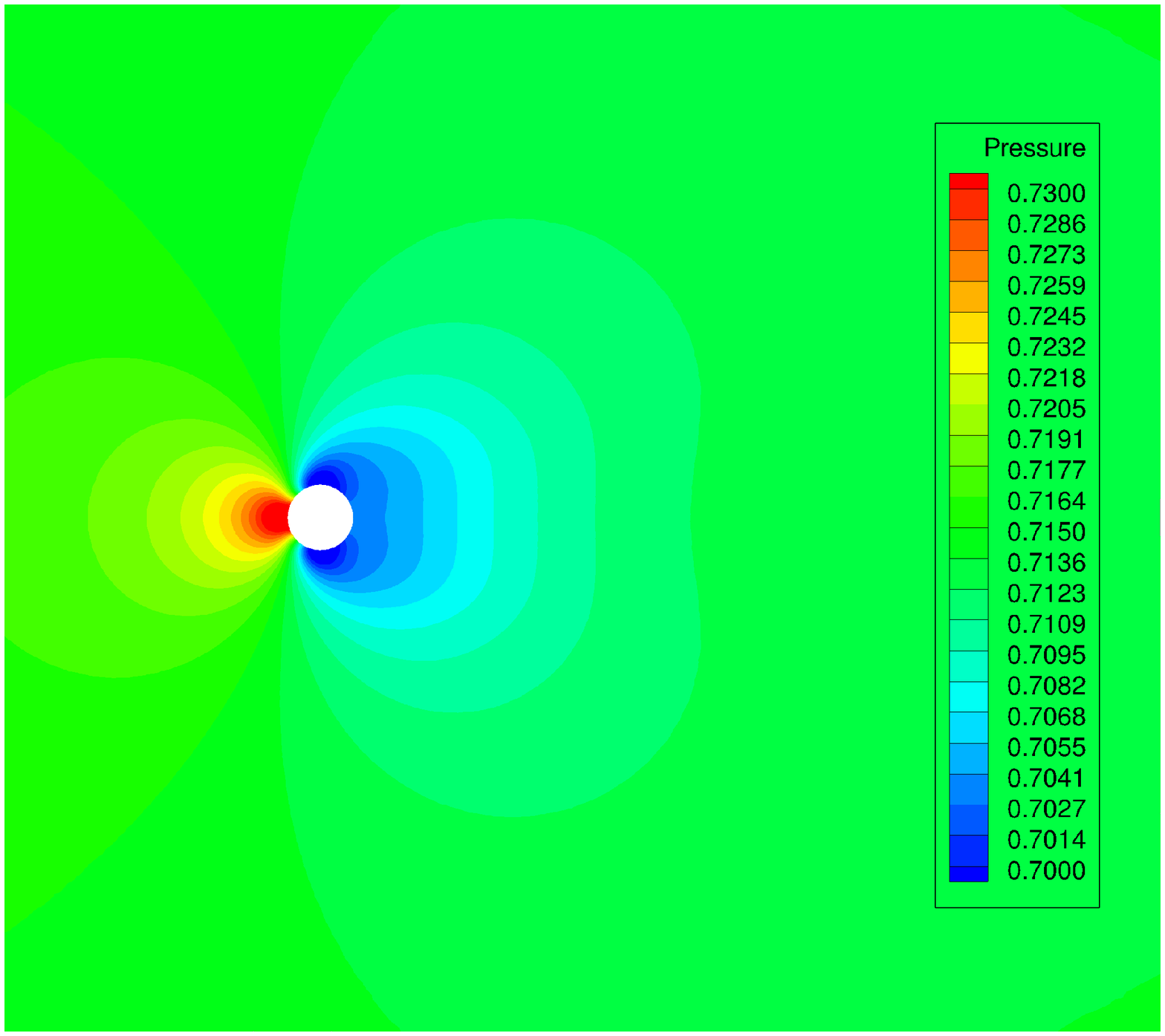}
		\caption{}
	\end{subfigure}
	\caption{The contour plots corresponding to the fifth order ($k=4$) solution in Example \ref{SubSec_SteadyCylinder}. 23 equally spaced contour levels are used to plot (a) Mach number between 0-0.22, and (b) pressure between 0.7-0.73}\label{Fig_cylinder}
\end{figure}

Following the notation in \cite{GAUTIER2013}, DDGIC results for $k=1,2,3,4$ are compared to selected data from the literature in Table \ref{Table_cylinder}. It should be noted that results presented in \citep{GAUTIER2013} correspond to incompressible flow simulations. This comparison is again justified by the fact that freestream Mach number is low, $M_{\infty}=0.2$, and thus, the flow behaves nearly incompressible. On other hand, the results in \citep{Canuto2015Tcvf} are obtained by compressible direct numerical simulations. It is seen that the second order ($k=1$) approximation overshoots the angle of separation while it undershoots $L_w/D$ and $a/D$. However, this situation is no longer the case for higher order approximations, and the wake characteristics obtained by the new DDGIC method for $k>1$ are in good agreement with the available data in the literature. For contour plots of Mach number and pressure corresponding to the fifth order ($k=4$) solution, see Figure \ref{Fig_cylinder}.

\begin{table}[h!]
\begin{center}
\caption{The comparison of the wake parameters between the new DDGIC method and the other results in the literature in Example \ref{SubSec_SteadyCylinder}}\label{Table_cylinder}
\begin{tabular}{cccccc}
\hlineB{3}
                                                & $C_d$     & $\theta$    & $L_w/D$   & $a/D$     & $b/D$   \\ \hlineB{2}
k=1                                             & 1.529     & 131.708     & 1.925     & 0.654     & 0.587     \\
k=2                                             & 1.567     & 127.499     & 2.246     & 0.711     & 0.593     \\
k=3                                             & 1.563     & 126.249     & 2.259     & 0.713     & 0.595     \\
k=4                                             & 1.564     & 126.151     & 2.258     & 0.713     & 0.595     \\ \hlineB{2}
Cheng et.al. \cite{cheng2018DDGIC}              & 1.561     & 126.200     & 2.180     & -         & -         \\
Canuto \& Daniel \cite{Canuto2015Tcvf}          & 1.560     & 126.300     & 2.280     & 0.724     & 0.596     \\
Gautier et.al.  \cite{GAUTIER2013}              & 1.48-1.62 & 124.4-127.3 & 2.13-2.35 & 0.71-0.76 & 0.59-0.60 \\ \hlineB{3}
\end{tabular}
\end{center}
\end{table}
 
\subsection{Unsteady flow over a cylinder}\label{SubSec_UnSteadyCylinder}

In this case, we use the same setup as in the previous example except for the Reynolds number $Re$ and the initial conditions. It is known that the flow over a cylinder becomes oscillatory after a critical Reynolds number of $Re_{cr}\approx 50$. For this purpose, we increased the Reynolds number to 75. The initial conditions are taken from the steady solutions of the previous example with a slightly perturbed y-velocity to ensure that unsteadiness occurs quickly.

\begin{figure}[!ht]
	\centering
	\includegraphics[scale=0.4]{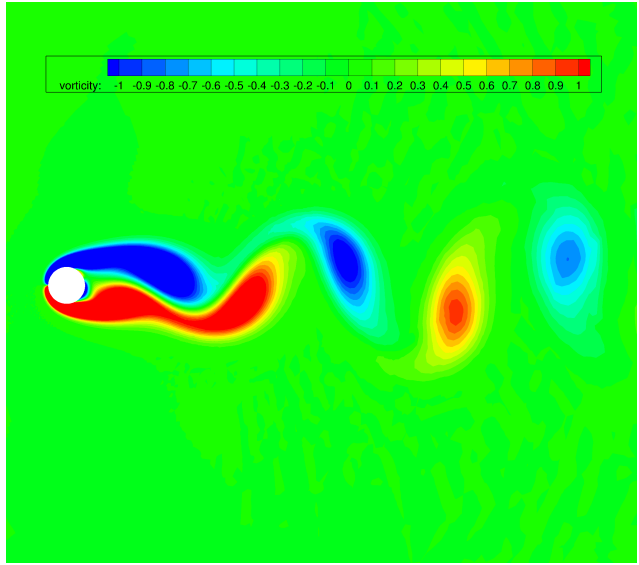}
	\caption{The contour plot corresponding to the fifth order ($k=4$) solution in Example \ref{SubSec_UnSteadyCylinder}. 21 equally spaced contour levels are used to plot nondimensional vorticity $\omega D/U_\infty$ between $\pm 1$}\label{Fig_cylinder_vortex}
\end{figure}

For this type of flows, one of the most popular flow features is the vortex shedding, which is also known as the Von Karman street in the literature. This phenomena is successfully captured by the new DDGIC method, and the corresponding nondimensional vorticity field for a fifth order solution ($k=4$) is depicted in Figure \ref{Fig_cylinder_vortex}.

As the vortices are shed in the wake region of the cylinder, the behavior of the flow becomes periodic. Therefore, the drag force exerted by the flow on the cylinder also becomes periodic. In addition, the obvious symmetry of the flow field in the steady example is no longer observed in this case. This induces a lift force on the cylinder, which is also oscillatory. Moreover, the periodicity of the flow allows us to compute the Strouhal number $St$, which is defined as

\begin{equation}
St=\frac{fD}{U}
\end{equation}
where $f$ is the frequency of the oscillations, $D$ is a characteristic length scale, i.e the diameter of the cylinder, and $U$ is a charactertistic velocity scale, i.e. the freestream velocity $U_\infty$. In Table \ref{Table_unsteady_cylinder}, the drag coefficient $C_d$,  the lift coefficient $C_l$, and the Strouhal number $St$ are reported. Note that the values for $C_d$ and $C_l$ correspond to the mean and deviations from the mean. It is seen that the second order ($k=1$) DDGIC method predicts slightly lower values for $C_d$, $C_l$, and $St$. On the other hand, the third and higher order solutions agree with each other as well as reported results in the litetarure.  

\begin{table}[h!]
\begin{center}
\caption{The comparison of the wake parameters between the new DDGIC method and the other results in the literature in Example \ref{SubSec_UnSteadyCylinder}}\label{Table_unsteady_cylinder}
\begin{tabular}{cccc}
\hlineB{3}
                                       & $C_d$              & $C_l$                    & $St$  \\ \hlineB{2}
$k=1$                                  & $1.309\pm 0.00261$ & $\pm 0.171$              & 0.146 \\
$k=2$                                  & $1.394\pm 0.00314$ & $\pm 0.209$              & 0.149 \\
$k=3$                                  & $1.396\pm 0.00342$ & $\pm 0.211$              & 0.150 \\
$k=4$                                  & $1.396\pm 0.00351$ & $\pm 0.211$              & 0.150 \\ \hlineB{2}
Canuto \& Daniel \cite{Canuto2015Tcvf} & $1.39\pm 0.00310$  & $\pm 0.212$ & 0.150 \\
Sun et al. \cite{SUN200641}           & -                  & -                        & 0.151 \\ \hlineB{3}
\end{tabular}
\end{center}
\end{table}

\section{Conclusions}\label{sec:Conclusion}

In this article, we have proposed a new formula to compute the nonlinear viscous numerical flux at the element edge and develop a new direct DG method with interface correction for compressible Navier-Stokes equations. The nonlinear viscous flux is decomposed into multiple individual diffusion processes and a new direction vector for each conserved variable is defined. The numerical flux formula of the original direct DG method is applied only to the conserved variable's gradient. The concept of multiple diffusion processes has been demonstrated to be useful for extending the new scheme to general equations and turbulence models. In a sequence of numerical examples, we have illustrated the optimal $(k+1)^{th}$ order convergence with high order $P_1$, $P_2$, $P_3$ and $P_4$ polynomial approximations. It has thus been shown that the new DDGIC method is capable of calculating the associated physical quantities accurately.

\newpage
\appendix
\section{The diffusion matrices for 2-D compressible Navier-Stokes equations}\label{sec:appendix-A}
In this section, we lay out the definitions of individual diffusion matrices $\bA^{(lm)}(\bQ)$ for 2-D compressible NS equations. Here, $\bA^{(lm)}(\bQ)$ should be understood as the diffusion matrix corresponding to the $l^{th}$ equation and the $m^{th}$ conserved variable.

According to the standard convention, continuity equation is the first equation, therefore we have $l=1$. Since there is no viscous terms in the continuity equation, we have that 
\begin{equation}
	\bA^{(11)}(\bQ)=\bA^{(12)}(\bQ)=\bA^{(13)}(\bQ)=\bA^{(14)}(\bQ)
		=\begin{pmatrix}
		  0 & 0 \\
		  0 & 0
		 \end{pmatrix}.
\end{equation} 
For the x-momentum equation, we have $l=2$ and
\begin{equation}
		\bA^{(21)}(\bQ) = 
		-\frac{\mu}{\rho}
		\begin{pmatrix}
  			\frac{4}{3}u & -\frac{2}{3}v \\
 			  v & u
		\end{pmatrix}, \,
		\bA^{(22)}(\bQ) =
		\frac{\mu}{\rho}
		\begin{pmatrix}
  			\frac{4}{3} & 0 \\
  			 0 & 1
		\end{pmatrix}, \,
		\bA^{(23)}(\bQ) =
		\frac{\mu}{\rho}
		\begin{pmatrix}
 			  0 & -\frac{2}{3} \\
 			  1 & 0
		\end{pmatrix},\,
		\bA^{(24)}(\bQ) =
		\begin{pmatrix}
 			  0 & 0\\
 			  0 & 0
		\end{pmatrix}. 
\end{equation}
Note that the total energy does not contribute to any of the diffusion matrices. In addition, one can easily show that all diffusion matrices for the x-momentum equation are compatible with 
\begin{equation}
	\bB^{(2)}(\bQ)=	
		\mu
		\begin{pmatrix}
  			\frac{4}{3}u & -\frac{2}{3}v \\
 			  v & u
		\end{pmatrix}.
\end{equation}
For the y-momentum equation, we have $l=3$ and
\begin{equation}
		\bA^{(31)}(\bQ) = 
		-\frac{\mu}{\rho}
		\begin{pmatrix}
  			v & u \\
 			-\frac{2}{3}u & \frac{4}{3}v
		\end{pmatrix}, \,
		\bA^{(32)}(\bQ) =
		\frac{\mu}{\rho}
		\begin{pmatrix}
 			  0 & 1 \\
 			  -\frac{2}{3} & 0
		\end{pmatrix}, \,
		\bA^{(33)}(\bQ) =
		\frac{\mu}{\rho}
		\begin{pmatrix}
  			 1 & 0 \\
  			 0 & \frac{4}{3}
		\end{pmatrix}, \,
		\bA^{(34)}(\bQ) =
		\begin{pmatrix}
 			  0 & 0\\
 			  0 & 0
		\end{pmatrix}. 
\end{equation}
As in the x-momentum equation, the total energy does not contribute to any of the diffusion matrices, and the diffusion matrices are all compatible with 
\begin{equation}
	\bB^{(3)}(\bQ)=	
		\mu
		\begin{pmatrix}
 			 v            & u   \\
  			-\frac{2}{3}u & \frac{4}{3}v \\
		\end{pmatrix}.
\end{equation}

For the total energy equation, we have $l=4$ and 
\begin{equation}
	\begin{aligned}
		\bA^{(41)}(\bQ) &= 
			\frac{\mu}{\rho}
			\begin{pmatrix}
				\left(\frac{\gamma}{2Pr}-\frac{4}{3}\right)u^2+ \left(\frac{\gamma}{2Pr}-1\right)v^2-\frac{\gamma}{Pr}e & -\frac{1}{3}uv \\ 
       	        -\frac{1}{3}uv  &	 \left(\frac{\gamma}{2Pr}-1\right)u^2+\left(\frac{\gamma}{2Pr}-\frac{4}{3}\right)v^2-\frac{\gamma}{Pr}e
			\end{pmatrix},  \\
		\bA^{(42)}(\bQ) &= 
			\frac{\mu}{\rho}
			\begin{pmatrix}
				\left(\frac{4}{3}-\frac{\gamma}{Pr}\right)u & v \\
				-\frac{2}{3}v & \left(1-\frac{\gamma}{Pr}\right)u
			\end{pmatrix},  \\
		\bA^{(43)}(\bQ) &= 
			\frac{\mu}{\rho}
			\begin{pmatrix}
				\left(1-\frac{\gamma}{Pr}\right)v &-\frac{2}{3}u \\
				u & \left(\frac{4}{3}-\frac{\gamma}{Pr}\right)v 
			\end{pmatrix},  \\
		\bA^{(44)}(\bQ)&=\frac{\mu\gamma}{\rho Pr}\mathbb{I}.
	\end{aligned}
\end{equation}
where $\mathbb{I}\in\mathbb{R}^{2\times 2}$ is the identity matrix. Unlike the momentum equations, $\bB^{(4)}(\bQ)$ does not exist. Therefore, the diffusion matrices of the total energy equation are not compatible.

\bibliographystyle{elsarticle-num-names}
\bibliography{Bib}

\begin{thebibliography}{49}
\providecommand{\natexlab}[1]{#1}
\providecommand{\url}[1]{\texttt{#1}}
\providecommand{\urlprefix}{URL }
\expandafter\ifx\csname urlstyle\endcsname\relax
  \providecommand{\doi}[1]{doi:\discretionary{}{}{}#1}\else
  \providecommand{\doi}[1]{doi:\discretionary{}{}{}\begingroup
  \urlstyle{rm}\url{#1}\endgroup}\fi
\providecommand{\bibinfo}[2]{#2}

\bibitem[{Liu and Yan(2010)}]{liuyan2010ddgic}
\bibinfo{author}{H.~Liu}, \bibinfo{author}{J.~Yan}, \bibinfo{title}{The direct
  discontinuous Galerkin ({DDG}) method for diffusion with interface
  corrections}, \bibinfo{journal}{Communications in Computational Physics}
  \bibinfo{volume}{8}~(\bibinfo{number}{3}) (\bibinfo{year}{2010})
  \bibinfo{pages}{541}.

\bibitem[{Lele(1992)}]{Lele1992}
\bibinfo{author}{S.~K. Lele}, \bibinfo{title}{Compact finite difference schemes
  with spectral-like resolution}, \bibinfo{journal}{Journal of Computational
  Physics} \bibinfo{volume}{103}~(\bibinfo{number}{1}) (\bibinfo{year}{1992})
  \bibinfo{pages}{16--42}.

\bibitem[{{Kopriva}(2009)}]{Kopriva-2009}
\bibinfo{author}{D.~A. {Kopriva}}, \bibinfo{title}{{Implementing spectral
  methods for partial differential equations. Algorithms for scientists and
  engineers}}, \bibinfo{publisher}{Berlin: Springer}, ISBN
  \bibinfo{isbn}{978-90-481-2260-8/hbk; 978-90-481-2261-5/ebook},
  \bibinfo{year}{2009}.

\bibitem[{Abgrall and Ricchiuto(2017)}]{Abgrall2017}
\bibinfo{author}{R.~Abgrall}, \bibinfo{author}{M.~Ricchiuto},
  \bibinfo{title}{High-Order Methods for CFD}, \bibinfo{publisher}{American
  Cancer Society}, ISBN \bibinfo{isbn}{9781119176817}, \bibinfo{pages}{1--54},
  \bibinfo{year}{2017}.

\bibitem[{Ekaterinaris(2005)}]{Ekaterinaris2005}
\bibinfo{author}{J.~A. Ekaterinaris}, \bibinfo{title}{High-order accurate, low
  numerical diffusion methods for aerodynamics}, \bibinfo{journal}{Progress in
  Aerospace Sciences} \bibinfo{volume}{41}~(\bibinfo{number}{3-4})
  (\bibinfo{year}{2005}) \bibinfo{pages}{192 -- 300}, ISSN
  \bibinfo{issn}{0376-0421}.

\bibitem[{Wang(2007)}]{Wang2007}
\bibinfo{author}{Z.~Wang}, \bibinfo{title}{High-order methods for the Euler and
  Navier-Stokes equations on unstructured grids}, \bibinfo{journal}{Progress in
  Aerospace Sciences} \bibinfo{volume}{43}~(\bibinfo{number}{1-3})
  (\bibinfo{year}{2007}) \bibinfo{pages}{1 -- 41}, ISSN
  \bibinfo{issn}{0376-0421}.

\bibitem[{Gassner(2009)}]{GasserThesis2009}
\bibinfo{author}{G.~Gassner}, \bibinfo{title}{Discontinuous Galerkin Methods
  for the Unsteady Compressible Navier Stokes Equations},
  \bibinfo{journal}{University of Stuttgart} \bibinfo{volume}{Dissertation}.

\bibitem[{Cockburn(2017)}]{Cockburn2017}
\bibinfo{author}{B.~Cockburn}, \bibinfo{title}{Discontinuous Galerkin Methods
  for Computational Fluid Dynamics}, \bibinfo{publisher}{American Cancer
  Society}, ISBN \bibinfo{isbn}{9781119176817}, \bibinfo{pages}{1--63},
  \bibinfo{year}{2017}.

\bibitem[{Wang et~al.(2013)Wang, Fidkowski, Abgrall, Bassi, Caraeni, Cary,
  Deconinck, Hartmann, Hillewaert, Huynh, Kroll, May, Persson, van Leer, and
  Visbal}]{Wang-Fidkowski-Abgrall-Bass-Caraeni-Cari-Deconinck-Hartmann2013}
\bibinfo{author}{Z.~Wang}, \bibinfo{author}{K.~Fidkowski},
  \bibinfo{author}{R.~Abgrall}, \bibinfo{author}{F.~Bassi},
  \bibinfo{author}{D.~Caraeni}, \bibinfo{author}{A.~Cary},
  \bibinfo{author}{H.~Deconinck}, \bibinfo{author}{R.~Hartmann},
  \bibinfo{author}{K.~Hillewaert}, \bibinfo{author}{H.~Huynh},
  \bibinfo{author}{N.~Kroll}, \bibinfo{author}{G.~May}, \bibinfo{author}{P.-O.
  Persson}, \bibinfo{author}{B.~van Leer}, \bibinfo{author}{M.~Visbal},
  \bibinfo{title}{High-order CFD methods: current status and perspective},
  \bibinfo{journal}{International Journal for Numerical Methods in Fluids}
  \bibinfo{volume}{72}~(\bibinfo{number}{8}) (\bibinfo{year}{2013})
  \bibinfo{pages}{811--845}.

\bibitem[{Cockburn et~al.(1998)Cockburn, Johnson, Shu, and
  Tadmor}]{Cockburn-book}
\bibinfo{author}{B.~Cockburn}, \bibinfo{author}{C.~Johnson},
  \bibinfo{author}{C.-W. Shu}, \bibinfo{author}{E.~Tadmor},
  \bibinfo{title}{Advanced numerical approximation of nonlinear hyperbolic
  equations}, vol. \bibinfo{volume}{1697} of \emph{\bibinfo{series}{Lecture
  Notes in Mathematics}}, \bibinfo{publisher}{Springer-Verlag},
  \bibinfo{address}{Berlin}, ISBN \bibinfo{isbn}{3-540-64977-8},
  \bibinfo{year}{1998}.

\bibitem[{Shu(2014)}]{Shu-DGreview2014}
\bibinfo{author}{C.-W. Shu}, \bibinfo{title}{Discontinuous {G}alerkin method
  for time-dependent problems: survey and recent developments}, in:
  \bibinfo{booktitle}{Recent developments in discontinuous {G}alerkin finite
  element methods for partial differential equations}, vol.
  \bibinfo{volume}{157} of \emph{\bibinfo{series}{IMA Vol. Math. Appl.}},
  \bibinfo{publisher}{Springer}, \bibinfo{pages}{25--62}, \bibinfo{year}{2014}.

\bibitem[{Hesthaven and Warburton(2007)}]{Hesthaven_nodal}
\bibinfo{author}{J.~S. Hesthaven}, \bibinfo{author}{T.~Warburton},
  \bibinfo{title}{Nodal Discontinuous Galerkin Methods: Algorithms, Analysis,
  and Applications}, \bibinfo{publisher}{Springer Publishing Company,
  Incorporated}, \bibinfo{edition}{1st} edn., ISBN \bibinfo{isbn}{0387720650},
  \bibinfo{year}{2007}.

\bibitem[{Rivi{\`e}re(2008)}]{Riviere2008}
\bibinfo{author}{B.~Rivi{\`e}re}, \bibinfo{title}{Discontinuous {G}alerkin
  methods for solving elliptic and parabolic equations},
  vol.~\bibinfo{volume}{35} of \emph{\bibinfo{series}{Frontiers in Applied
  Mathematics}}, \bibinfo{publisher}{Society for Industrial and Applied
  Mathematics (SIAM), Philadelphia, PA}, \bibinfo{year}{2008}.

\bibitem[{Di~Pietro and Ern(2012)}]{Pietro-Ern2011-DGbook}
\bibinfo{author}{D.~A. Di~Pietro}, \bibinfo{author}{A.~Ern},
  \bibinfo{title}{Mathematical aspects of discontinuous {G}alerkin methods},
  vol.~\bibinfo{volume}{69} of \emph{\bibinfo{series}{Math\'ematiques \&
  Applications (Berlin) [Mathematics \& Applications]}},
  \bibinfo{publisher}{Springer, Heidelberg}, \bibinfo{year}{2012}.

\bibitem[{Zhang and Shu(2011)}]{zhang-MPS-Review}
\bibinfo{author}{X.~Zhang}, \bibinfo{author}{C.-W. Shu},
  \bibinfo{title}{Maximum-principle-satisfying and positivity-preserving
  high-order schemes for conservation laws: survey and new developments},
  \bibinfo{journal}{Proc. R. Soc. A} ~(\bibinfo{number}{467})
  (\bibinfo{year}{2011}) \bibinfo{pages}{2752--2776.}

\bibitem[{Bassi and Rebay(1997)}]{BR1}
\bibinfo{author}{F.~Bassi}, \bibinfo{author}{S.~Rebay}, \bibinfo{title}{A
  High-Order Accurate Discontinuous Finite Element Method for the Numerical
  Solution of the Compressible Navier–Stokes Equations},
  \bibinfo{journal}{Journal of Computational Physics}
  \bibinfo{volume}{131}~(\bibinfo{number}{2}) (\bibinfo{year}{1997})
  \bibinfo{pages}{267--279}.

\bibitem[{Bassi et~al.(2005)Bassi, Crivellini, Rebay, and Savini}]{BR4}
\bibinfo{author}{F.~Bassi}, \bibinfo{author}{A.~Crivellini},
  \bibinfo{author}{S.~Rebay}, \bibinfo{author}{M.~Savini},
  \bibinfo{title}{Discontinuous Galerkin solution of the Reynolds-averaged
  Navier--Stokes and k--$\omega$ turbulence model equations},
  \bibinfo{journal}{Computers \& Fluids}
  \bibinfo{volume}{34}~(\bibinfo{number}{4-5}) (\bibinfo{year}{2005})
  \bibinfo{pages}{507--540}.

\bibitem[{Nguyen et~al.(2011)Nguyen, Peraire, and
  Cockburn}]{Nguyen-Peraire-Cockburn2011}
\bibinfo{author}{N.~C. Nguyen}, \bibinfo{author}{J.~Peraire},
  \bibinfo{author}{B.~Cockburn}, \bibinfo{title}{High-order implicit
  hybridizable discontinuous {G}alerkin methods for acoustics and
  elastodynamics}, \bibinfo{journal}{J. Comput. Phys.}
  \bibinfo{volume}{230}~(\bibinfo{number}{10}) (\bibinfo{year}{2011})
  \bibinfo{pages}{3695--3718}.

\bibitem[{Peraire and Persson(2008)}]{peraire2008compact}
\bibinfo{author}{J.~Peraire}, \bibinfo{author}{P.-O. Persson},
  \bibinfo{title}{The compact discontinuous Galerkin (CDG) method for elliptic
  problems}, \bibinfo{journal}{SIAM Journal on Scientific Computing}
  \bibinfo{volume}{30}~(\bibinfo{number}{4}) (\bibinfo{year}{2008})
  \bibinfo{pages}{1806--1824}.

\bibitem[{Persson(2013)}]{Persson2013}
\bibinfo{author}{P.-O. Persson}, \bibinfo{title}{A sparse and high-order
  accurate line-based discontinuous {G}alerkin method for unstructured meshes},
  \bibinfo{journal}{J. Comput. Phys.} \bibinfo{volume}{233}
  (\bibinfo{year}{2013}) \bibinfo{pages}{414--429}.

\bibitem[{Zhang(2017)}]{Zhang-2017-NS}
\bibinfo{author}{X.~Zhang}, \bibinfo{title}{On positivity-preserving high order
  discontinuous {G}alerkin schemes for compressible {N}avier-{S}tokes
  equations}, \bibinfo{journal}{J. Comput. Phys.} \bibinfo{volume}{328}
  (\bibinfo{year}{2017}) \bibinfo{pages}{301--343}.

\bibitem[{Arnold(1982)}]{arnold1982interior}
\bibinfo{author}{D.~N. Arnold}, \bibinfo{title}{An interior penalty finite
  element method with discontinuous elements}, \bibinfo{journal}{SIAM journal
  on numerical analysis} \bibinfo{volume}{19}~(\bibinfo{number}{4})
  (\bibinfo{year}{1982}) \bibinfo{pages}{742--760}.

\bibitem[{Wheeler(1978)}]{wheeler1978elliptic}
\bibinfo{author}{M.~F. Wheeler}, \bibinfo{title}{An elliptic collocation-finite
  element method with interior penalties}, \bibinfo{journal}{SIAM Journal on
  Numerical Analysis} \bibinfo{volume}{15}~(\bibinfo{number}{1})
  (\bibinfo{year}{1978}) \bibinfo{pages}{152--161}.

\bibitem[{Baker(1977)}]{baker1977finite}
\bibinfo{author}{G.~A. Baker}, \bibinfo{title}{Finite element methods for
  elliptic equations using nonconforming elements},
  \bibinfo{journal}{Mathematics of Computation}
  \bibinfo{volume}{31}~(\bibinfo{number}{137}) (\bibinfo{year}{1977})
  \bibinfo{pages}{45--59}.

\bibitem[{Hartmann and Houston(2006{\natexlab{a}})}]{hartmann2006a}
\bibinfo{author}{R.~Hartmann}, \bibinfo{author}{P.~Houston},
  \bibinfo{title}{Symmetric interior penalty DG methods for the compressible
  Navier-Stokes equations I: Method Formulation},
  \bibinfo{journal}{International Journal of Numerical Analysis \& Modeling}
  \bibinfo{volume}{3}~(\bibinfo{number}{1})
  (\bibinfo{year}{2006}{\natexlab{a}}) \bibinfo{pages}{1--20}.

\bibitem[{Hartmann and Houston(2006{\natexlab{b}})}]{hartmann2006b}
\bibinfo{author}{R.~Hartmann}, \bibinfo{author}{P.~Houston},
  \bibinfo{title}{Symmetric interior penalty DG methods for the compressible
  Navier--Stokes equations II: Goal--oriented a posteriori error estimation},
  \bibinfo{journal}{International Journal of Numerical Analysis \& Modeling}
  \bibinfo{volume}{3}~(\bibinfo{number}{2})
  (\bibinfo{year}{2006}{\natexlab{b}}) \bibinfo{pages}{141--162}.

\bibitem[{Hartmann and Houston(2008)}]{hartmann2008}
\bibinfo{author}{R.~Hartmann}, \bibinfo{author}{P.~Houston}, \bibinfo{title}{An
  optimal order interior penalty discontinuous Galerkin discretization of the
  compressible Navier--Stokes equations}, \bibinfo{journal}{Journal of
  Computational Physics} \bibinfo{volume}{227}~(\bibinfo{number}{22})
  (\bibinfo{year}{2008}) \bibinfo{pages}{9670--9685}.

\bibitem[{Liu and Yan(2008)}]{liuyan2008ddg}
\bibinfo{author}{H.~Liu}, \bibinfo{author}{J.~Yan}, \bibinfo{title}{THE DIRECT
  DISCONTINUOUS GALERKIN ({DDG}) METHODS FOR DIFFUSION PROBLEMS},
  \bibinfo{journal}{SIAM Journal on Numerical Analysis}
  \bibinfo{volume}{47}~(\bibinfo{number}{1}) (\bibinfo{year}{2008})
  \bibinfo{pages}{675--698}.

\bibitem[{Vidden and Yan(2013)}]{vidden2013sddg}
\bibinfo{author}{C.~Vidden}, \bibinfo{author}{J.~Yan}, \bibinfo{title}{A NEW
  DIRECT DISCONTINUOUS GALERKIN METHOD WITH SYMMETRIC STRUCTURE FOR NONLINEAR
  DIFFUSION EQUATIONS}, \bibinfo{journal}{Journal of Computational Mathematics}
  \bibinfo{volume}{31}~(\bibinfo{number}{6}) (\bibinfo{year}{2013})
  \bibinfo{pages}{638--662}.

\bibitem[{Yan(2013)}]{yan2013new}
\bibinfo{author}{J.~Yan}, \bibinfo{title}{A new nonsymmetric discontinuous
  Galerkin method for time dependent convection diffusion equations},
  \bibinfo{journal}{Journal of Scientific Computing}
  \bibinfo{volume}{54}~(\bibinfo{number}{2}) (\bibinfo{year}{2013})
  \bibinfo{pages}{663--683}.

\bibitem[{Chen et~al.(2016)Chen, Huang, and Yan}]{chen2016}
\bibinfo{author}{Z.~Chen}, \bibinfo{author}{H.~Huang},
  \bibinfo{author}{J.~Yan}, \bibinfo{title}{Third order
  maximum-principle-satisfying direct discontinuous Galerkin methods for time
  dependent convection diffusion equations on unstructured triangular meshes},
  \bibinfo{journal}{Journal of Computational Physics} \bibinfo{volume}{308}
  (\bibinfo{year}{2016}) \bibinfo{pages}{198 -- 217}, ISSN
  \bibinfo{issn}{0021-9991}.

\bibitem[{Zhang and Yan(2017)}]{Zhang-Yan-2017}
\bibinfo{author}{M.~Zhang}, \bibinfo{author}{J.~Yan}, \bibinfo{title}{Fourier
  type super convergence study on {DDGIC} and symmetric {DDG} methods},
  \bibinfo{journal}{J. Sci. Comput.}
  \bibinfo{volume}{73}~(\bibinfo{number}{2-3}) (\bibinfo{year}{2017})
  \bibinfo{pages}{1276--1289}.

\bibitem[{Huang et~al.(2020)Huang, Li, and Yan}]{Huang-Yan-2020}
\bibinfo{author}{H.~Huang}, \bibinfo{author}{J.~Li}, \bibinfo{author}{J.~Yan},
  \bibinfo{title}{High order symmetric direct discontinuous {G}alerkin method
  for elliptic interface problems with fitted mesh}, \bibinfo{journal}{J.
  Comput. Phys.} \bibinfo{volume}{409} (\bibinfo{year}{2020})
  \bibinfo{pages}{109301, 23}.

\bibitem[{Kannan and Wang(2010)}]{KANNAN20102007}
\bibinfo{author}{R.~Kannan}, \bibinfo{author}{Z.~Wang}, \bibinfo{title}{The
  direct discontinuous Galerkin (DDG) viscous flux scheme for the high order
  spectral volume method}, \bibinfo{journal}{Computers \& Fluids}
  \bibinfo{volume}{39}~(\bibinfo{number}{10}) (\bibinfo{year}{2010})
  \bibinfo{pages}{2007 -- 2021}.

\bibitem[{Yang et~al.(2016)Yang, Cheng, Wang, Luo, Si, and
  Pandare}]{yang2016fast}
\bibinfo{author}{X.~Yang}, \bibinfo{author}{J.~Cheng},
  \bibinfo{author}{C.~Wang}, \bibinfo{author}{H.~Luo}, \bibinfo{author}{J.~Si},
  \bibinfo{author}{A.~Pandare}, \bibinfo{title}{A fast, implicit discontinuous
  Galerkin method based on analytical Jacobians for the compressible
  Navier-Stokes equations}, in: \bibinfo{booktitle}{54th AIAA Aerospace
  Sciences Meeting}, \bibinfo{pages}{1326}, \bibinfo{year}{2016}.

\bibitem[{Cheng et~al.(2016{\natexlab{a}})Cheng, Yang, Liu, Liu, and
  Luo}]{cheng2016DDG}
\bibinfo{author}{J.~Cheng}, \bibinfo{author}{X.~Yang},
  \bibinfo{author}{X.~Liu}, \bibinfo{author}{T.~Liu}, \bibinfo{author}{H.~Luo},
  \bibinfo{title}{A direct discontinuous Galerkin method for the compressible
  Navier–Stokes equations on arbitrary grids}, \bibinfo{journal}{Journal of
  Computational Physics} \bibinfo{volume}{327}
  (\bibinfo{year}{2016}{\natexlab{a}}) \bibinfo{pages}{484 -- 502}.

\bibitem[{Cheng et~al.(2016{\natexlab{b}})Cheng, Liu, Yang, Liu, and
  Luo}]{cheng2016DDG_RANS}
\bibinfo{author}{J.~Cheng}, \bibinfo{author}{X.~Liu},
  \bibinfo{author}{X.~Yang}, \bibinfo{author}{T.~Liu},
  \bibinfo{author}{H.~Luo}, \bibinfo{title}{A Direct Discontinuous Galerkin
  Method for Computation of Turbulent Flows on Hybrid Grids}
  \doi{\bibinfo{doi}{10.2514/6.2016-3333}}.

\bibitem[{Cheng et~al.(2017)Cheng, Liu, Liu, and Luo}]{cheng_liu_liu_luo_2017}
\bibinfo{author}{J.~Cheng}, \bibinfo{author}{X.~Liu}, \bibinfo{author}{T.~Liu},
  \bibinfo{author}{H.~Luo}, \bibinfo{title}{A Parallel, High-Order Direct
  Discontinuous Galerkin Method for the Navier-Stokes Equations on 3D Hybrid
  Grids}, \bibinfo{journal}{Communications in Computational Physics}
  \bibinfo{volume}{21}~(\bibinfo{number}{5}) (\bibinfo{year}{2017})
  \bibinfo{pages}{1231–1257}.

\bibitem[{Yue et~al.(2017)Yue, Cheng, and Liu}]{yue_cheng_liu_2017}
\bibinfo{author}{H.~Yue}, \bibinfo{author}{J.~Cheng}, \bibinfo{author}{T.~Liu},
  \bibinfo{title}{A Symmetric Direct Discontinuous Galerkin Method for the
  Compressible Navier-Stokes Equations}, \bibinfo{journal}{Communications in
  Computational Physics} \bibinfo{volume}{22}~(\bibinfo{number}{2})
  (\bibinfo{year}{2017}) \bibinfo{pages}{375–392}.

\bibitem[{Cheng et~al.(2018{\natexlab{a}})Cheng, Yue, Yu, and
  Liu}]{cheng2018DDGIC}
\bibinfo{author}{J.~Cheng}, \bibinfo{author}{H.~Yue}, \bibinfo{author}{S.~Yu},
  \bibinfo{author}{T.~Liu}, \bibinfo{title}{A Direct Discontinuous Galerkin
  Method with Interface Correction for the Compressible Navier-Stokes Equations
  on Unstructured Grids}, \bibinfo{journal}{Advances in applied mathematics and
  mechanics} \bibinfo{volume}{10}~(\bibinfo{number}{1})
  (\bibinfo{year}{2018}{\natexlab{a}}) \bibinfo{pages}{1--21}.

\bibitem[{Yang et~al.(2018)Yang, Cheng, Luo, and Zhao}]{YANG2018216}
\bibinfo{author}{X.~Yang}, \bibinfo{author}{J.~Cheng},
  \bibinfo{author}{H.~Luo}, \bibinfo{author}{Q.~Zhao}, \bibinfo{title}{A
  reconstructed direct discontinuous Galerkin method for simulating the
  compressible laminar and turbulent flows on hybrid grids},
  \bibinfo{journal}{Computers \& Fluids} \bibinfo{volume}{168}
  (\bibinfo{year}{2018}) \bibinfo{pages}{216 -- 231}.

\bibitem[{Cheng et~al.(2018{\natexlab{b}})Cheng, Yue, Yu, and
  Liu}]{CHENG2018adaptive}
\bibinfo{author}{J.~Cheng}, \bibinfo{author}{H.~Yue}, \bibinfo{author}{S.~Yu},
  \bibinfo{author}{T.~Liu}, \bibinfo{title}{Analysis and development of
  adjoint-based h-adaptive direct discontinuous Galerkin method for the
  compressible Navier–Stokes equations}, \bibinfo{journal}{Journal of
  Computational Physics} \bibinfo{volume}{362}
  (\bibinfo{year}{2018}{\natexlab{b}}) \bibinfo{pages}{305 -- 326}.

\bibitem[{Xiaoquan et~al.(2019)Xiaoquan, Cheng, Luo, and Zhao}]{Xiaoquan2019}
\bibinfo{author}{Y.~Xiaoquan}, \bibinfo{author}{J.~Cheng},
  \bibinfo{author}{H.~Luo}, \bibinfo{author}{Q.~Zhao}, \bibinfo{title}{Robust
  Implicit Direct Discontinuous Galerkin Method for Simulating the Compressible
  Turbulent Flows}, \bibinfo{journal}{AIAA Journal}
  \bibinfo{volume}{57}~(\bibinfo{number}{3}) (\bibinfo{year}{2019})
  \bibinfo{pages}{1113--1132}.

\bibitem[{Gottlieb et~al.(2001)Gottlieb, Shu, and
  Tadmor}]{Gottlieb-Shu-Tabmor-2001}
\bibinfo{author}{S.~Gottlieb}, \bibinfo{author}{C.-W. Shu},
  \bibinfo{author}{E.~Tadmor}, \bibinfo{title}{Strong stability-preserving
  high-order time discretization methods}, \bibinfo{journal}{SIAM Rev.}
  \bibinfo{volume}{43}~(\bibinfo{number}{1}) (\bibinfo{year}{2001})
  \bibinfo{pages}{89--112}.

\bibitem[{Shu and Osher(1988)}]{Shu-Osher-1988}
\bibinfo{author}{C.-W. Shu}, \bibinfo{author}{S.~Osher},
  \bibinfo{title}{Efficient implementation of essentially non-oscillatory
  shock-capturing schemes}, \bibinfo{journal}{Journal of Computational Physics}
  \bibinfo{volume}{77}~(\bibinfo{number}{2}) (\bibinfo{year}{1988})
  \bibinfo{pages}{439 -- 471}.

\bibitem[{Mengaldo et~al.(2014)Mengaldo, Grazia, Witherden, Farrington,
  Vincent, Sherwin, and Peiro}]{nektar_code}
\bibinfo{author}{G.~Mengaldo}, \bibinfo{author}{D.~D. Grazia},
  \bibinfo{author}{F.~Witherden}, \bibinfo{author}{A.~Farrington},
  \bibinfo{author}{P.~Vincent}, \bibinfo{author}{S.~Sherwin},
  \bibinfo{author}{J.~Peiro}, \bibinfo{title}{A Guide to the Implementation of
  Boundary Conditions in Compact High-Order Methods for Compressible
  Aerodynamics}, \bibinfo{publisher}{7th AIAA Theoretical Fluid Mechanics
  Conference}, \bibinfo{year}{2014}.

\bibitem[{Gautier et~al.(2013)Gautier, Biau, and Lamballais}]{GAUTIER2013}
\bibinfo{author}{R.~Gautier}, \bibinfo{author}{D.~Biau},
  \bibinfo{author}{E.~Lamballais}, \bibinfo{title}{A reference solution of the
  flow over a circular cylinder at Re=40}, \bibinfo{journal}{Computers \&
  Fluids} \bibinfo{volume}{75} (\bibinfo{year}{2013}) \bibinfo{pages}{103 --
  111}, ISSN \bibinfo{issn}{0045-7930}.

\bibitem[{Canuto and Taira(2015)}]{Canuto2015Tcvf}
\bibinfo{author}{D.~Canuto}, \bibinfo{author}{K.~Taira},
  \bibinfo{title}{Two-dimensional compressible viscous flow around a circular
  cylinder}, \bibinfo{journal}{Journal of fluid mechanics}
  \bibinfo{volume}{785} (\bibinfo{year}{2015}) \bibinfo{pages}{349--371},
  \doi{\bibinfo{doi}{10.1017/jfm.2015.635}}.

\bibitem[{Sun et~al.(2006)Sun, Wang, and Liu}]{SUN200641}
\bibinfo{author}{Y.~Sun}, \bibinfo{author}{Z.~Wang}, \bibinfo{author}{Y.~Liu},
  \bibinfo{title}{Spectral (finite) volume method for conservation laws on
  unstructured grids VI: Extension to viscous flow}, \bibinfo{journal}{Journal
  of Computational Physics} \bibinfo{volume}{215}~(\bibinfo{number}{1})
  (\bibinfo{year}{2006}) \bibinfo{pages}{41 -- 58}, ISSN
  \bibinfo{issn}{0021-9991}.

\end{thebibliography}

\end{document}